\documentclass[onefignum,onetabnum]{siamart220329}
\usepackage{graphics,graphicx,epstopdf,url}
\usepackage{arydshln}
\usepackage{color}
\usepackage{geometry}
\usepackage{algorithmic}
\usepackage{url,cases,supertabular}
\usepackage{amssymb,mathtools}
\usepackage{enumerate,makecell}
\usepackage{amsfonts}
\usepackage{graphicx,epstopdf}
\usepackage{color,verbatim}
\usepackage{pifont}
\usepackage{mathrsfs,subeqnarray,subfigure}
\usepackage{listings,array}
\usepackage{booktabs}
\usepackage{dashrule}
\usepackage{arydshln}
\usepackage{graphicx}
\usepackage{amsfonts}
\usepackage{mathrsfs}
\usepackage{graphicx}  \usepackage{epstopdf}
\usepackage{subfigure}  \usepackage{lineno}
\usepackage{multirow,bm}
\usepackage{color}
\usepackage{algorithm}
\usepackage{algorithmic}
\usepackage{longtable}
\usepackage{mdwtab}
\usepackage{mdwmath}
\usepackage{bbding}
\usepackage{booktabs}
\usepackage{rotating}
\usepackage{enumitem}
\usepackage{lineno}
\usepackage{url}
\usepackage{pdflscape}
\usepackage{makecell}

\newcommand{\revise}[1]{#1}
\newcommand{\abstractrevision}[1]{{\color{blue}#1}}
\newenvironment{revisionblock}{}{}

\numberwithin{equation}{section}

\newtheorem{assumption}[theorem]{Assumption}

\newcommand{\ba}{\begin{array}}
\newcommand{\ea}{\end{array}}

\newcommand{\bit}{\begin{itemize}}
\newcommand{\eit}{\end{itemize}}
\newcommand{\be}{\begin{equation}}
\newcommand{\ee}{\end{equation}}
\newcommand{\bee}{\begin{equation*}}
\newcommand{\eee}{\end{equation*}}
\newcommand{\bea}{\begin{eqnarray}}
\newcommand{\eea}{\end{eqnarray}}

\newcommand{\st}{\mathrm{s.t.}}

\newcommand{\by}{\mathrm{\bf y}}
\newcommand{\bs}{\mathrm{\bf s}}
\newcommand{\bx}{\mathbf{x}}
\newcommand{\bd}{\mathbf{d}}
\newcommand{\bq}{\mathbf{q}}

\newcommand{\bv}{\mathbf{v}}
\newcommand{\bW}{\mathbf{W}}

\newcommand{\Rmn}[1]{\uppercase\expandafter{\romannumeral#1}}
\renewcommand{\theenumi}{\roman{enumi}}

\newcommand{\Pcal}{\mathcal{P}}

\numberwithin{equation}{section}

\newcommand{\Mcal}{\mathcal{M}}

\newcommand{\grad}{\mathrm{grad}}

\newcommand{\R}{\mathbb{R}}

\numberwithin{theorem}{section}
\newcommand{\iprod}[2]{\left \langle #1, #2 \right \rangle }

\usepackage{lipsum}
\usepackage{clrscode}
\usepackage{appendix}
\usepackage{bm}
\usepackage{booktabs}
\usepackage{url}
\usepackage{multirow}
\usepackage{textcomp}
\usepackage{amsmath}
\usepackage{amsfonts}
\usepackage{amssymb}
\usepackage{mathrsfs}
\usepackage{hyperref}
\usepackage{graphicx,graphics,subfigure}
\usepackage{hyperref,url}
\usepackage{epsf,epstopdf}
\usepackage{algorithm}
\usepackage{algorithmic}
\usepackage{bbm}
\usepackage{dsfont}
\usepackage{indentfirst}
\usepackage{pdfpages}
\usepackage{pdflscape}
\usepackage{longtable}
\ifpdf
  \DeclareGraphicsExtensions{.eps,.pdf,.png,.jpg}
\else
  \DeclareGraphicsExtensions{.eps}
\fi
\allowdisplaybreaks[1]

\headers{Decentralized projected Riemannian gradient method}{Kangkang Deng, Jiang Hu}

\title{Decentralized projected Riemannian stochastic recursive momentum method for smooth optimization on compact submanifolds}

\author{Kangkang Deng\thanks{ Department of Mathematics,  National University of Defense Technology, Changsha, 410073,
China (\email{freedeng1208@gmail.com}).}
\and Jiang Hu\thanks{Corresponding author. Department of Mathematics, University of California, Berkeley, CA 94720, US
(\email{hujiangopt@gmail.com}).}
}

\usepackage{amsopn}

\ifpdf
\hypersetup{
  pdftitle={Decentralized projected Riemannian gradient method},
  pdfauthor={Kangkang Deng, Jiang Hu}
}
\fi

\allowdisplaybreaks[2]
\begin{document}

\maketitle

\begin{abstract}
This paper studies decentralized optimization over a compact submanifold within a communication network of $n$ nodes, where each node possesses a smooth non-convex local cost function, and the goal is to jointly minimize the sum of these local costs. We focus particularly on the online setting, where local data is processed in real-time as it streams in, without the need for full data storage. We propose a decentralized projected Riemannian stochastic recursive momentum (DPRSRM) method that employs local hybrid stochastic gradient estimators and uses the network to track the global gradient. DPRSRM achieves an oracle complexity of  \(\mathcal{O}(\epsilon^{-\frac{3}{2}})\), outperforming existing methods that have at most \(\mathcal{O}(\epsilon^{-2})\) complexity. Our method requires only $\mathcal{O}(1)$ gradient evaluations per iteration for each local node and does not require restarting with a large batch gradient.   Furthermore, we demonstrate the effectiveness of our proposed methods compared to state-of-the-art ones through numerical experiments on principal component analysis problems and low-rank matrix completion. \abstractrevision{This revised version corrects the mean-square oracle condition and the estimator/tracking analysis, and replaces the original separate error summation by a coupled Lyapunov argument; Algorithm~\ref{alg:drgta} and the $\mathcal O(\epsilon^{-3/2})$ oracle-complexity order are unchanged.}

\end{abstract}


\section{Introduction}

Decentralized optimization has emerged as a prominent area of research, particularly for its application in large-scale systems, such as sensor networks, distributed computing, and machine learning. In these contexts, data is often partitioned across numerous nodes, rendering centralized optimization approaches impractical due to challenges such as privacy limitations and restricted computational resources. In this work, we are concerned with the distributed smooth optimization on a compact submanifold
\be \label{prob:original}
\begin{aligned}
  \min \quad &  \frac{1}{n}\sum_{i=1}^n f_i(x_i), \\
  \st \quad & x_1 = \cdots = x_n, \;\;\\
  & x_i \in \Mcal, \;\; \forall~ i=1,2,\ldots, n,
\end{aligned}
\ee
where $n$ is the number of nodes, $f_i$ is the smooth nonconvex local objective at the $i$-th node, and $\Mcal$ is a (nonconvex) compact smooth submanifold {embedded in} $\R^{d\times r}$ {with the extrinsic dimensions $(d,r)$}, e.g., the Stiefel manifold ${\rm St}(d,r):=\{ x \in \R^{d\times r} : x^\top x = I_r \}$.  
Problem \eqref{prob:original} is prevalent in machine learning, signal processing, and deep learning, see, e.g., the principal component analysis \cite{ye2021deepca}, the low-rank matrix completion \cite{boumal2015low,kasai2019riemannian,deng2023decentralized}, the low-dimensional subspace learning \cite{ando2005framework,mishra2019riemannian}, and the deep neural networks with batch normalization \cite{cho2017riemannian,hu2022riemannian}. One challenge in solving \eqref{prob:original} comes from the nonconvexity of the manifold constraint, which causes difficulty in achieving the consensus \cite{chen2021decentralized,deng2023decentralized,hu2023decentralized,deng2023decentralized11}.

In this paper, we investigate an online setting where each node \(i\) interacts with its local cost function \(f_i\) through a stochastic first-order oracle (SFO).  This SFO setting is particularly relevant in various online learning and expected risk minimization problems, where the noise introduced by the SFO stems from the variability of sampling over streaming data received at each node. A notable example is online principal component analysis \cite{cardot2018online}. Our primary focus is on the oracle complexity, defined as the total number of SFO queries required at each node to compute an \(\epsilon\)-stationary point tuple \(\{x_1, \cdots, x_n\}\), as formalized in Definition \ref{def:station}.


\subsection{Related works}


Decentralized optimization in Euclidean space (i.e., $\mathcal{M} = \mathbb{R}^{d \times r}$) has been extensively studied over the past few decades (see, e.g., \cite{bianchi2012convergence,shi2015extra,xu2015augmented,qu2017harnessing,di2016next,tatarenko2017non, wai2017decentralized,yuan2018exact,hong2017prox,zeng2018nonconvex, scutari2019distributed, sun2020improving}). However, since problem \eqref{prob:original} involves a manifold $\Mcal$, which is often nonlinear and nonconvex, these works may fail when directly applied to solve problem \eqref{prob:original}. 



Perhaps the earliest works for solving \eqref{prob:original} are  \cite{mishra2019riemannian,shah2017distributed}. However, these methods require an asymptotically infinite number of consensus steps for convergence, which limits their practical applicability.
For the case where $\Mcal$ is the Stiefel manifold, \cite{chen2021decentralized} propose a decentralized Riemannian gradient descent method and its gradient-tracking version.   To use a single step of consensus, augmented Lagrangian methods \cite{wang2022decentralized,wang2022variance} are also investigated, where a different stationarity is used.  \cite{sun2024global} propose a decentralized retraction-free gradient tracking algorithm, and show that it exhibits ergodic
$\mathcal{O}(1/K)$ convergence rate.   
However, these studies rely on the orthogonal structure of the Stiefel manifold.  Recently, \cite{deng2023decentralized} used the projection operators instead of retractions and expanded the distributed Riemannian gradient descent algorithm and the gradient tracking version to the compact submanifolds of Euclidean space. 
Moreover, the integration of decentralized manifold optimization with other algorithms has also been proposed, including the conjugate gradient algorithm \cite{chen2023decentralized} and the natural gradient method \cite{hu2023decentralized}. 
Furthermore, \cite{hu2024improving} achieves single-step consensus for the general compact submanifold by carefully elaborating on the smoothness structure and the asymptotic 1-Lipschitz continuity of the projection operator associated with the submanifold geometry.

Several studies have focused on the finite-sum setting of problem \eqref{prob:original}, where $f_i = \frac{1}{m}\sum_{r = 1}^{m_i} f_{i,r}$.  
\cite{chen2021decentralized} propose a decentralized  Riemannian stochastic gradient descent method. By combining the variable sample size gradient approximation method with the gradient tracking dynamic, 
\cite{zhao2024distributed} propose a distributed Riemannian stochastic optimization algorithm on the Stiefel manifold. Although both methods can also be used in the online setting,  the oracle complexities are $\mathcal{O}(\epsilon^{-2})$, which is not optimal. 
 It is worth noting that the decentralized variance reduced method \cite{wang2022variance} has been studied. However, they need to periodically calculate the full gradient, which is not suitable for the online setting.


\subsection{Contribution} 
In this paper, we propose DPRSRM, a novel online variance-reduced method for decentralized non-convex manifold optimization with stochastic first-order oracles (SFO). 

\begin{itemize}

   \item To achieve
fast and robust performance, the DPRSRM algorithm is built upon gradient tracking \cite{chen2021decentralized,deng2023decentralized} and a stochastic gradient momentum estimator \cite{cutkosky2019momentum,han2020riemannian}, which can be viewed as online variance reduction method. The only existing decentralized stochastic variance-reduced manifold optimization algorithm is the VRSGT proposed by \cite{wang2022variance}. Note that VRSGT is a double-loop algorithm that requires very large minibatch sizes. Conversely, the proposed DPRSRM is a single-loop algorithm with $\mathcal{O}(1)$ oracle queries per update. Numerical experiments demonstrate the effectiveness of the proposed methods compared to state-of-the-art ones through eigenvalue problems and low-rank matrix completion.

    \item  Our algorithm achieves an oracle complexity of $\mathcal{O}(\epsilon^{-\frac{3}{2}})$.   A comparison of the oracle complexity of DPRSRM with related algorithms is provided in Table \ref{tab:my_label}, where DPRSRM achieves a lower oracle complexity than the existing decentralized stochastic manifold optimization algorithms \cite{deng2023decentralized,chen2021decentralized,wang2022variance,zhao2024distributed}. 
    Moreover, DPRSRM uses a single step of consensus to achieve communication, compared to other project/retraction algorithms \cite{chen2021decentralized,deng2023decentralized,zhao2024distributed} that need $\log_{\sigma_2}(\frac{1}{2\sqrt{n}})$ rounds of consensus, where $\sigma_2$ is the second largest singular value of the communication graph matrix.  
\end{itemize}

\begin{table*}
    \centering
    \begin{tabular}{|c|c|c|c|c|c|}
  \hline
        Algorithm  & Manifold types & Communication & Tracking & VR  & Oracle \\ \hline
           \cite{chen2021decentralized}     &   Stiefel manifold &  multiple  & \ding{55} & \ding{55} & $\mathcal{O}(\epsilon^{-2})$  \\ \hline
 \cite{zhao2024distributed}     &  Stiefel manifold & multiple  & \ding{55} & 
 \ding{55} & $\mathcal{O}(\epsilon^{-2})$\\\hline
         This paper  &  compact submanifold & single  & \checkmark & \checkmark & $\mathcal{O}(\epsilon^{-3/2})$\\ \hline
    \end{tabular}
    \caption{Comparison of the oracle complexity results of Riemannian online decentralized methods.  ``Communication'' means rounds of communications per iteration.  ``Tracking'' denotes the gradient tracking, ``VR'' denotes variance reduction.  We do not list the work in \cite{wang2022variance} since they focus on the finite-sum setting and are not applicable to the online setting. }
    \label{tab:my_label}
\end{table*}

\subsection{Notation}
For the compact submanifold $\Mcal$ of $\mathbb{R}^{d \times r}$, we always take the Euclidean metric $\iprod{\cdot}{\cdot}$ as the Riemannian metric. We use $\|\cdot \|$ to denote the Euclidean norm. We denote the $n$-fold Cartesian product of $\Mcal$ as $\Mcal^n = \Mcal \times \cdots \times \Mcal$.
For any $x\in \Mcal$, the tangent space and normal space of $\Mcal$ at $x$ are denoted by $T_x\Mcal$ and $N_x\Mcal$, respectively. For a differentiable function $h: \R^{d\times r} \rightarrow \R$, we denote its Euclidean gradient by $\nabla h(x)$ and its Riemannian gradient by $\grad h(x)$.  
For a positive integer $n$, define $[n] =\{1, \dots, n\}$. Let $\mathbf{1}_n\in \mathbb{R}^n$ be a vector where all entries are equal to 1. Define $J := \frac{1}{n}\mathbf{1}_n\mathbf{1}_n^{\top}$.   Unless otherwise explicitly defined, we now provide explanations for all lowercase variables used in this paper. Take $x$ as an example, we denote $x_i$ as a local variable at $i$-th node; $\hat{x} = \frac{1}{n}\sum_{i=1}^nx_i$ is the Euclidean average. Moreover, we use the bold notations $\bx := [x_1^\top, \dots, x_n^\top]^\top \in \R^{nd \times r},~~\hat{\bx} := [\hat{x}^\top, \dots, \hat{x}^\top]^\top \in \R^{nd \times r},$ where $\bx$ denotes the collection of all local variables $x_i$ and $\hat{\bx}$ is $n$ copies of $\hat{x}$. When applied to the iterative process, in $k$-th iteration, we use $x_{i,k}$ to denote a local variable at $i$-th node and $\hat{x}_k = \frac{1}{n}\sum_{i=1}^nx_{k,i}$. Similarly, we also 
denote 
$\bx_k := [x_{1,k}^\top, \dots, x_{n,k}^\top]^\top \in \R^{nd \times r},~~~ \hat{\bx}_k = [\hat{x}_k^\top, \dots, \hat{x}_k^\top]^\top \in \R^{nd \times r}.$ 
Other lowercase variables can also be denoted similarly as $x$. \revise{Whenever the projection is single-valued, we write $\bar{x}:=\Pcal_{\Mcal}(\hat{x})$ and $\bar{\bx}:=\mathbf 1_n\otimes\bar{x}$.} Define the function $f(\bx): = \sum_{i=1}^n f_i(x_i)$. \revise{For a stacked argument, $\grad f(\bx)$ denotes $[\grad f_1(x_1)^\top,\ldots,\grad f_n(x_n)^\top]^\top$; for a single point $z\in\Mcal$, $f(z):=\frac1n\sum_{i=1}^n f_i(z)$.} Let $\bW := W \otimes I_d \in \R^{nd \times nd}$, where $\otimes$ denotes the Kronecker product. 
        

\section{Preliminary}
This section introduces the definition of stationary point for problem \eqref{prob:original}, and gives a key property for compact submanifolds, i.e., proximal smoothness. 
\subsection{Stationary point}
Let $x_1,\cdots,x_n\in \Mcal$ represent the local copies of each node. 
Let $\Pcal_{\Mcal}$ be the orthogonal projection onto $\Mcal$. Note that for $\{x_i\}_{i=1}^n \subset \Mcal$,
\[ {\rm argmin}_{y\in\Mcal}\sum_{i=1}^n \|y - x_i\|^2 = \Pcal_{\Mcal}(\hat{x}). \]

Any element $\bar{x}$ in $\Pcal_{\Mcal}(\hat{x})$ is the induced arithmetic mean of $\{x_i\}_{i=1}^n$ on $\Mcal$ \cite{sarlette2009consensus}.
\revise{Thus, for a single point $z\in\Mcal$, $\grad f(z)=\frac1n\sum_{i=1}^n\grad f_i(z)$.} The $\epsilon$-stationary point of problem \eqref{prob:original} is defined as follows.
\begin{definition}[\cite{deng2023decentralized}] \label{def:station}
The set of points $\{x_1,x_2,\cdots,x_n\} \subset \Mcal$ is called an $\epsilon$-stationary point of \eqref{prob:original} if there exists a $\bar{x} \in \Pcal_{\Mcal}(\hat{x})$ such that
\[\mathbb{E}\left[ \frac{1}{n}\sum_{i=1}^n \| x_i - \bar{x}\|^2 \right] \leq \epsilon \quad {\rm and} \quad \mathbb{E}[\|\grad f(\bar{x})\|^2] \leq \epsilon. \]
\end{definition}
We refer to these two terms as consensus error and optimality error, respectively. 

\subsection{Proximal smoothness}\label{sec:proximal}
Proximal smoothness is an effective tool for addressing the nonconvex nature of manifold constraints within decentralized optimization settings \cite{deng2023decentralized}. Define the distance from a point \(x \in \mathbb{R}^{d \times r}\) to the manifold \(\mathcal{M}\) and the nearest-point projection of \(x\) onto \(\mathcal{M}\) as $
\text{dist}(x, \mathcal{M}) := \inf_{y \in \mathcal{M}} \|y - x\| \;\; \text{and} \;\; \mathcal{P}_{\mathcal{M}}(x) := \arg\min_{y \in \mathcal{M}} \|y - x\|,$
respectively. For a given number \(R > 0\), the \(R\)-tube around \(\mathcal{M}\) is defined as the set $U_{\mathcal{M}}(R) := \{x : \text{dist}(x, \mathcal{M}) < R\}.$

A closed set \(\mathcal{M}\) is said to be \(R\)-proximally smooth if the projection \(\mathcal{P}_{\mathcal{M}}(x)\) is unique whenever \(\text{dist}(x, \mathcal{M}) < R\). Following \cite{clarke1995proximal},  an $R$-proximally smooth set $\mathcal{M}$ satisfies that for any real ${\delta}\in (0,R)$,
    \be \label{proximally-1}
    \left\| \Pcal_{\mathcal{M}} (x) -\Pcal_{\mathcal{M}} (y)\right\| \leq \frac{R}{R-{\delta}}\|x - y\|,~~ \forall x,y \in \bar{U}_{\mathcal{M}}({\delta}),
    \ee
where \(\bar{U}_{\mathcal{M}}(\delta) := \{x : \text{dist}(x, \mathcal{M}) \leq \delta\}\). For instance, the Stiefel manifold is known to be a 1-proximally smooth set \cite{balashov2021gradient}, {which also provides additional examples of rank manifolds along with their specific proximal smoothness radii.}

\section{Problem setup and the proposed DPRSRM}\label{sec:main-alg}
In this section, we present the problem setup considered in this paper, outlining the assumptions for the objective function and the communication network. Building on this setup, we then develop a decentralized algorithm to solve problem \eqref{prob:original} and provide the main convergence rate results. 
\subsection{Problem setup}
Let us start with {assumptions on problem \eqref{prob:original}, based on the fact that any compact $C^2$-submanifolds in Euclidean space belong to proximally smooth set \cite{balashov2021gradient,clarke1995proximal,davis2020stochastic}.}  
\begin{assumption}\label{assum-f}
    {Assume that $\Mcal$ is proximally smooth with radius $R$.} Each objective function $f_i$ is gradient Lipschitz continuous with modulus $L_{f}$ on the convex hull of $\Mcal$, denoted by ${\rm conv}(\Mcal)$. \revise{Moreover, the aggregate objective $f(x):=\frac1n\sum_{i=1}^n f_i(x)$ is bounded below on $\Mcal$; denote its minimum by $f_*$.}  
\end{assumption}
Under this assumption, the following Riemannian quadratic upper bound for $f_i$ has been estabilshed in \cite{deng2023decentralized}. 
\begin{lemma}[\cite{deng2023decentralized}, Lemma 2]\label{lemma:lipsctz}
Under Assumption \ref{assum-f}, there exists $L_g$, for any $x,y\in \Mcal$, the following inequality holds:
\begin{equation}\label{f-riemannian-Lip}
    f_i(y) - f_i(x) \leq \left<\grad f_i(x),y-x\right>+\frac{L_g}{2}\|y-x\|^2, ~ i\in [n].
\end{equation}   
Moreover, there exists a constant $L_G>0$ such that
\begin{equation}\label{g-riemannian-lip}
    \| \grad f_i(x) - \grad f_i(y) \| \leq L_G \|x - y\|,~i\in [n].
\end{equation}
\end{lemma}

We now present the assumption for the communication network. 
Denote by the undirected graph $G:=\{ \mathcal{V}, \mathcal{E}\}$, where $\mathcal{V}=\{1,2,\ldots, n\}$ is the set of all nodes and $\mathcal{E}$ is the set of edges. {Let $W$ be the mixing matrix associated with $G$ and $W_{ij} \ne 0$ if there is an edge $(i,j) \in \mathcal{E}$ and $W_{ij} = 0$ otherwise.} We use the following standard assumptions on $W$, see, e.g., \cite{chen2021decentralized,zeng2018nonconvex}.
\begin{assumption} \label{assum-w}
    We assume that the undirected graph $G$ is connected and $W$ is doubly stochastic, i.e., (i) $W=W^{\top}$; (ii) $W_{i j} \geq 0$ and $1>W_{i i}>0$ for all $i, j$; (iii) the eigenvalues of $W$ lie in $(-1,1]$. \revise{In addition, the second-largest singular value $\sigma_2$ of $W$ satisfies $\sigma_2\in[0,1)$.}
\end{assumption}
Consider a sufficiently rich probability space \(\{\Omega, \mathbb{P}, \mathcal{F}\}\). For a given decentralized algorithm, we assume it generates an iterative sequence \(\{x_{i,k}\}_{k \geq 0}\), where \(x_{i,k}\) denotes the \(k\)-th iteration at node \(i\). At each step, node \(i\) observes a random vector \(\xi_{i,k}\). We then define an increasing sequence of sub-\(\sigma\)-algebras within \(\mathcal{F}\), constructed from the random vectors observed in succession by the network nodes: for all $k \geq 1$,
$
\mathcal{F}_0 := \{\Omega, \emptyset\}, \quad \mathcal{F}_k := \sigma\left(\{\xi_{i,0}, \xi_{i,1}, \ldots, \xi_{i,k-1} : i \in \mathcal{V}\}\right),
$
where \(\emptyset\) represents the empty set. The following assumptions are made regarding the stochastic gradient \(\nabla f_i(x, \xi_{i,k})\):
\begin{assumption}\label{assum:stochatis}
\begin{revisionblock}
For any $\mathcal{F}_k$--measurable variable $x\in\Mcal$ and $k\geq 0$, the algorithm generates a sample $\xi_{i,k}\sim \Omega$ for each node $i$ and returns a stochastic gradient $\grad f_i(x,\xi_{i,k})$. There are constants $\nu_i\geq0$ and $\bar L>0$; we write $\nu^2:=\frac{1}{n}\sum_{i=1}^n\nu_i^2$. The oracle satisfies
\begin{align}
  &  \mathbb{E}\left[ \grad f_i(x,\xi_{i,k}) \vert \mathcal{F}_k \right] = \grad f_i(x), \label{assm:mean} \\ 
 &    \mathbb{E}\left[ \|\grad f_i(x,\xi_{i,k}) - \grad f_i(x) \|^2 \vert \mathcal{F}_k\right] \leq \nu_i^2. \label{assm:vari}
\end{align}
Moreover, the collection $\{\xi_{i,k}:i\in \mathcal{V},k\geq0\}$ consists of independent random variables and the same-sample stochastic gradients are mean-square $\bar L$-Lipschitz:
\begin{equation}\label{assm:ms-lipschitz}
   \mathbb{E}\!\left[ \| \grad f_i(x,\xi_{i,k}) - \grad f_i(y,\xi_{i,k}) \|^2 \mid \mathcal{F}_k\right] \leq \bar{L}^2\|x-y\|^2.
\end{equation}
Finally, the un-clipped initialization used in Algorithm~\ref{alg:drgta} satisfies $\|\grad f_i(x_{i,0},\xi_{i,0})\|\leq B$ almost surely.
\end{revisionblock}
\end{assumption}

To measure the oracle complexity,  we give the definition of a stochastic first-order oracle (SFO) for \eqref{prob:original}.
\begin{definition}[\textbf{stochastic first-order oracle}]\label{def:stochastic}
    For the problem \eqref{prob:original}, a stochastic first-order oracle for each node $i$ is defined as follows: compute the Riemannian gradient $\grad f_i(x,\xi_i)$ given a sample $\xi_i\in \Omega$.  
\end{definition}

\subsection{The Algorithm}

In this subsection, we introduce the Decentralized Projected Riemannian Stochastic Recursive Momentum (DPRSRM) method for addressing \eqref{prob:original}, along with its associated convergence results. Inspired by the notable effectiveness of variance reduction and gradient tracking in decentralized frameworks \cite{wang2022variance,deng2023decentralized,chen2021decentralized}, we seek a novel combination of variance reduction and gradient tracking for decentralized online problems on compact submanifolds for improving the oracle complexity. In particular, we focus on the following stochastic gradient estimator using the momentum technique introduced in \cite{cutkosky2019momentum,han2020riemannian,deng2024oracle}:
\begin{equation}\label{update:qik}
\begin{aligned}
    q_{i,k} = &   \grad f_{i}(x_{i,k},\xi_{i,k}) \\
   & +(1-\tau)(d_{i,k-1} - \grad f_{i}(x_{i,k-1},\xi_{i,k})).
\end{aligned}
\end{equation}

We note that \eqref{update:qik} can be rewritten as 
\begin{equation}\label{update:qik1}
\begin{aligned}
    q_{i,k} = &  \tau \grad f_{i}(x_{i,k},\xi_{i,k})  +(1-\tau)(d_{i,k-1} \\
    &+\grad f_{i}(x_{i,k},\xi_{i,k})- \grad f_{i}(x_{i,k-1},\xi_{i,k})),
\end{aligned}
\end{equation}
which hybrids stochastic gradient $\grad f_{i}(x_{i,k},\xi_{i,k})$ with the recursive gradient estimator in RSARAH/SRG \cite{han2020variance,zhou2019faster} for $\tau \in[0,1]$.  Since the direction $q_{i,k}$ may be unbounded, we introduce a clipped gradient estimator $d_{i,k}$
\begin{equation}\label{update:dik}
    d_{i,k} = \left\{
    \begin{array}{cc}
        q_{i,k} & \text{if} ~ \| q_{i,k} \| \leq B, \\
        B \frac{q_{i,k}}{\|q_{i,k}\|} & \text{otherwise},    \end{array}
    \right.
\end{equation}
where $B>0$ is a user-defined constant.   To further reduce the variance of the gradient estimator in different nodes, we compute the gradient tracking iteration as follows:
\be\label{update:sik}
s_{i,k}  = \sum_{j=1}^n W_{ij}s_{j,k-1} + d_{i,k} - d_{i,k-1}, ~i\in [n].
\ee

 A crucial advantage of gradient tracking-type methods lies in the applicability of the use of a constant step size $\alpha$. Since $s_{i,k}$ may not remain in the tangent space $T_{x_{i,k}}\Mcal$, we introduce {its tangent space projection} $v_{i,k} = \Pcal_{T_{x_{i,k}} (s_{i,k})}$ and update the new iterate $x_{i,k+1}$ as follows:
 \begin{equation}\label{updata-new-}
     x_{i,k+1} = \Pcal_{\Mcal}(\sum_{j=1}^n W_{ij}x_{j,k} - \alpha
v_{i,k} ), ~i\in [n],
 \end{equation}
 {where the projection-based average \cite{deng2023decentralized,hu2024improving}, i.e., $\Pcal_{\Mcal}(\sum_{j=1}^n W_{ij}x_{j,k})$, is adopted to ensure the decrease in the consensus error among nodes.}
For ease of analysis, we stack variables in each node and rewrite \eqref{update:qik}, \eqref{update:sik} and \eqref{updata-new-} as 
\begin{equation}
\left\{
    \begin{aligned}
        \bq_k  =& \grad f(\bx_k,\mathbf{\xi}_k) \\
        &+ (1-\tau) (\bd_{k-1} - \grad f(\bx_{k-1},\mathbf{\xi}_k)) \\
        \bs_{k} =& \bW\bs_{k-1} + \bd_k - \bd_{k-1} \\
        \bv_k  = &\Pcal_{T_{\bx}\Mcal^n}( \bs_k ) \\
        \bx_{k+1}  =&  \Pcal_{\Mcal^n}(\bW \bx_k - \alpha \bv_k),
    \end{aligned}
    \right.
\end{equation}
where  $\Pcal_{T_{\bx}\Mcal^n}:= \Pcal_{T_{x_{1,k}}\Mcal} \times \dots \times \Pcal_{T_{x_{n,k}}\Mcal}$ and $\mathbf{\xi}_k = \{\xi_{i,k}\}_{i\in \mathcal{V}}$. The overall algorithm is given in Algorithm \ref{alg:drgta}.

\begin{algorithm}[htbp]
\caption{The DPRSRM for solving \eqref{prob:original}} \label{alg:drgta}
\begin{algorithmic}[1]
\REQUIRE  \revise{Initial point $\bx_{0} = \bar{\bx}_0 \in \Mcal^n$, $\bs_{-1} =  \bd_{-1} = \mathbf{0}$, $\alpha,\tau$, and clipping radius $B$.}
\STATE Sample $\xi_{i,0}$,  let $s_{i,0} = d_{i,0} = \grad f_i(x_{i,0},\xi_{i,0})$.
\STATE $v_{i,0} = \Pcal_{T_{x_{i,0}}\Mcal}( s_{i,0})$.
\STATE $ x_{i,1} = \Pcal_{\Mcal}(\sum_{j=1}^n W_{ij}x_{j,0} - \alpha
v_{i,0} )$.
\FOR{$k = 1,2,\cdots,K$}
\STATE Update stochastic gradient estimator $q_{i,k}$ via \eqref{update:qik}
\STATE Update the clipped gradient estimator $d_{i,k}$ via \eqref{update:dik}.
\STATE Update Riemannian gradient tracking $s_{i,k}$ via \eqref{update:sik}.
\STATE Project onto tangent space: $v_{i,k} = \Pcal_{T_{x_{i,k}}\Mcal}( s_{i,k})$.
\STATE Update new iterate $x_{i,k+1}$ via \eqref{updata-new-}. 
\ENDFOR
\end{algorithmic}
\end{algorithm}



\subsection{Main results}\label{sec:main-resilt}

The convergence analysis of Algorithm \ref{alg:drgta} can be divided into two parts: the consensus error and the optimality error, as defined in Definition  \ref{def:station}.  Let us first focus on the consensus error.  We define a neighborhood around $\bx\in \Mcal^n$ as follows:
\begin{equation}\label{def:N-neiborhood-1}
    \mathcal{N}(\delta):=\{ \bx\in \Mcal^n: \|\bx - \bar{\bx} \| \leq \delta \}.
\end{equation}
Note that if $\bx\in \mathcal{N}(\delta)$, then $x_{i}\in \bar{U}_{\Mcal}(\delta)$ for any $i\in [n]$. 

\revise{By appropriately selecting the step size $\alpha$ and initializing with $\bx_0 \in \mathcal{N}(\delta)$, we can establish an explicit relationship between the consensus error and the step size. The proof is given in Section \ref{sec:linear-consensus}.}
\begin{revisionblock}
Define
\begin{equation}\label{def:repaired-constants}
 L:=\max\{L_g,L_G,\bar L\},\qquad
 G_{\max}:=\max_{i\in[n],\,x\in\Mcal}\|\grad f_i(x)\|,
\end{equation}
where $G_{\max}<\infty$ follows from compactness. Throughout the convergence
analysis we take $B\geq G_{\max}$ and set
\[
 D_B:=\left(\frac{3B}{1-\sigma_2}\right)^2,
 \qquad \kappa_M:=\frac{R}{R-2\delta},
 \qquad \rho:=\kappa_M\sigma_2.
\]
\begin{theorem}[Consensus error] \label{lem:consensus}
Let $\{\bx_k\}_k$ be generated by Algorithm~\ref{alg:drgta} and suppose that
Assumptions~\ref{assum-f}--\ref{assum:stochatis} hold. Let $B\geq G_{\max}$
and suppose that
   \begin{equation}\label{cond-alpha-delta}
\begin{aligned}
    \alpha & \leq  \min\left\{ \frac{\delta}{\sqrt{nD_B}},
      \frac{(R(1-\sigma_2)-2\delta)\delta}{R\sqrt{nD_B}} \right\}, \\
    \delta &<  \frac{1}{2}\min\left\{ R,  R(1-\sigma_2)\right\}.
    \end{aligned}
\end{equation}
If $\bx_0=\bar\bx_0$, then every projection in Algorithm~\ref{alg:drgta} is
single-valued in the prescribed tube and, for all $k\geq0$,
    \begin{equation}\label{eq:xk-barxk-bound}
    \frac{1}{n} \|\bar{\bx}_k - \bx_k \|^2 \leq C_B\alpha^2,
    \qquad C_B:=\frac{\kappa_M^2D_B}{(1-\rho)^2}.
    \end{equation}
\end{theorem}
   
We now state the repaired optimality result. Its proof is given in
Section~\ref{sec:proof:dprgt}.
\begin{theorem}[Optimality error] \label{thm:dprgt}
Let the hypotheses of Theorem~\ref{lem:consensus} hold, with $\alpha$ chosen
below. There exists a
problem-dependent constant $\kappa_0>0$. Choose a fixed
$\kappa\geq\kappa_0$ and then a sufficiently small problem-dependent constant
$\theta>0$. There are constants $K_0,C_{\rm opt}>0$, independent of the
horizon $K$, such that, for every $K\geq K_0$, setting
\begin{equation}\label{eq:repaired-parameters}
 \alpha=\theta(K+1)^{-1/3},\qquad \tau=\kappa\alpha^2\leq1,
\end{equation}
with $\alpha$ also satisfying \eqref{cond-alpha-delta}. Then
\begin{equation}\label{eq:repaired-optimality-rate}
 \min_{0\leq k\leq K}\mathbb{E}\|\grad f(\bar x_k)\|^2
 \leq C_{\rm opt}(K+1)^{-2/3}.
\end{equation}
\end{theorem}

\begin{corollary}\label{coll}
Under the conditions of Theorem~\ref{thm:dprgt}, there is an index
$k_*\leq K$ such that
\[
 \mathbb{E}\|\grad f(\bar x_{k_*})\|^2=\mathcal O(K^{-2/3}),
 \qquad
 \frac1n\mathbb{E}\|\bx_{k_*}-\bar\bx_{k_*}\|^2
 =\mathcal O(K^{-2/3}).
\]
Since each node evaluates only a constant number of stochastic gradients per
iteration, DPRSRM has per-node SFO complexity
$\mathcal O(\epsilon^{-3/2})$ for Definition~\ref{def:station}.
\end{corollary}
\end{revisionblock}

According to Corollary \ref{coll}, DPRSRM achieves an oracle complexity of  $\mathcal{O}(\epsilon^{-\frac{3}{2}})$, outperforming existing methods \cite{chen2021decentralized,zhao2024distributed}, which have an oracle complexity of at most \(\mathcal{O}(\epsilon^{-2})\).


\section{Outline of convergence analysis}\label{sec4}
As shown in Section \ref{sec:main-resilt},  the convergence analysis consists of two critical components: the consensus error and the optimality measure, corresponding to Theorem \ref{lem:consensus} and Theorem \ref{thm:dprgt}, respectively. In this section, we will outline the proof. 
 \subsection{Consensus error}  \label{sec:linear-consensus}

This subsection addresses the consensus error in Theorem \ref{lem:consensus}. 
In Algorithm \ref{alg:drgta}, the update of the main iterate $x_{k+1}$ involves the projection on $\Mcal$. Due to the nonconvex nature of compact submanifolds, the projection is not always unique. Before proving Theorem \ref{lem:consensus}, we will first demonstrate that the projection operator in Algorithm \ref{alg:drgta} is well-defined, meaning that the points being projected are always ensured to be within a neighborhood that belongs to  $\bar{U}_{\Mcal}(R)$.  

\begin{revisionblock}
We first investigate the uniform boundedness of $\|\bs_k\|$ in the following lemma.
\begin{lemma}\label{lemma:neibohood}
Let $\{\bx_k\}_k$ be generated by Algorithm~\ref{alg:drgta}. Under
Assumptions~\ref{assum-w}--\ref{assum:stochatis}, for every $k\geq0$,
\be\label{eq:sk-bound}
 \|\bs_k\|^2 \leq nD_B.
\ee
\end{lemma}

 The following lemma shows that the iterates $\{\bx_k\}$ remain in
 $\mathcal{N}(\delta)$ under the stated conditions.

\begin{lemma} \label{lem:stay-neighborhood}

Suppose that Assumptions~\ref{assum-f}--\ref{assum:stochatis} hold and
$B\geq G_{\max}$. Let $x_{i,k}$ be generated by Algorithm~\ref{alg:drgta},
and let $\alpha$ and $\delta$ satisfy \eqref{cond-alpha-delta}.
If $\bx_0=\bar\bx_0$, then for every $k\geq0$,
$\bx_k\in\mathcal N(\delta)$ and
\begin{align}
    \sum_{j=1}^n {W_{ij}} x_{j,k} -\alpha v_{i,k} & \in \bar{U}_{\Mcal}(2\delta), ~ i=1,\cdots,n.
\end{align}
\end{lemma}

Now we give the proof of Theorem \ref{lem:consensus}. 

\begin{proof}[Proof of Theorem \ref{lem:consensus}]
Since $\bx_0=\bar\bx_0$, Lemma~\ref{lem:stay-neighborhood} gives, for every
$k\geq0$, $\bx_k\in\mathcal N(\delta)$ and
\begin{equation}\label{eq:neibohood3}
    \begin{aligned}
       &  \sum_{j=1}^n {W_{ij}} x_{j,k} -\alpha v_{i,k}  \in \bar{U}_{\Mcal} (2\delta), ~ i = 1,\cdots,n.
    \end{aligned}
\end{equation}
By the definition of $\bar{\bx}_{k+1}$, we have
\[ \begin{aligned}
    \|\bx_{k+1} - \bar{\bx}_{k+1}\| & \leq \| \bx_{k+1} - \bar{\bx}_k \| \\
    & = \| \Pcal_{\Mcal^n}(\bW \bx_k - \alpha  \bv_k) - \Pcal_{\Mcal^n}( \hat{\bx}_k) \| \\
    & \leq \frac{R}{R-2\delta} \|\bW \bx_k - \alpha  \bv_k -\hat{\bx}_k  \| \\
    & \leq \rho\|\bx_k-\bar{\bx}_k\|
      +\kappa_M\sqrt{nD_B}\,\alpha,
\end{aligned}
\]
where the first inequality follows from the optimality of $\bar{\bx}_{k+1}$,
and the second uses \eqref{eq:neibohood3} and the
$\kappa_M$-Lipschitz continuity of $\Pcal_{\Mcal}$ over
$\bar U_{\Mcal}(2\delta)$. Hence
 \be \label{eq:temp:pgd:consensus}
    \begin{aligned}
        \|\bar{\bx}_{k+1}-\bx_{k+1}\|
        &\leq \rho^{k+1}\|\bar{\bx}_0-\bx_0\|
        +\frac{\kappa_M\sqrt{nD_B}}{1-\rho}\alpha.
    \end{aligned}
    \ee   
The initialization $\bx_0=\bar\bx_0$ proves \eqref{eq:xk-barxk-bound}.
\end{proof}
\end{revisionblock}

\subsection{Optimality error}\label{sec:proof:dprgt}

\begin{revisionblock}
We replace the separate summation of the estimation and tracking errors by a
single coupled argument. Define the normalized deterministic sequences
\begin{equation}\label{eq:five-errors}
\begin{aligned}
 \mathsf{G}_k&:=\mathbb{E}\|\grad f(\bar x_k)\|^2,&
 \mathsf{S}_k&:=\frac1n\mathbb{E}\|\bs_k\|^2,\\
 \mathsf{E}_k&:=\frac1n\mathbb{E}
   \|\bd_k-\grad f(\bx_k)\|^2,&
 \mathsf{T}_k&:=\frac1n\mathbb{E}\|\bs_k-\hat{\bd}_k\|^2,\\
 \mathsf{X}_k&:=\frac1n\mathbb{E}\|\bx_k-\bar\bx_k\|^2.&
\end{aligned}
\end{equation}

\begin{lemma}[Coupled one-step inequalities]\label{lem:coupled-one-step}
Suppose that Assumptions~\ref{assum-f}--\ref{assum:stochatis} hold,
$B\geq G_{\max}$, $0<\tau\leq1$, and the local projection conditions in
Theorem~\ref{lem:consensus} are satisfied. There are constants $d_X,C_D>0$,
independent of $k,K,\alpha,\tau$, such that
\begin{equation}\label{eq:coupled-E-main}
 \mathsf{E}_{k+1}\le(1-\tau)^2\mathsf{E}_k+64\bar L^2\mathsf{X}_k
      +16\bar L^2\alpha^2\mathsf{S}_k+2\tau^2\nu^2.
\end{equation}
\begin{equation}\label{eq:coupled-T-main}
 \mathsf{T}_{k+1}\le q_T\mathsf{T}_k+h_T(96\bar L^2\mathsf{X}_k
      +24\bar L^2\alpha^2\mathsf{S}_k
      +3\tau^2\mathsf{E}_k+3\tau^2\nu^2).
\end{equation}
\begin{equation}\label{eq:coupled-X-main}
 \mathsf{X}_{k+1}\le q_X\mathsf{X}_k+c_X\alpha^2\mathsf{S}_k.
\end{equation}
\begin{equation}\label{eq:coupled-D-main}
 \begin{aligned}
 \mathbb{E} f(\bar x_{k+1})
 &\le\mathbb{E} f(\bar x_k)-\frac\alpha4\mathsf{G}_k
      -\frac\alpha4\mathsf{S}_k
      +\frac{3\alpha}{2}(\mathsf{E}_k+\mathsf{T}_k)\\
 &\quad+d_X\alpha\mathsf{X}_k+C_D\alpha^3.
 \end{aligned}
\end{equation}
where
\[
 q_T:=\frac{1+\sigma_2^2}{2},\quad
 h_T:=\frac{1+\sigma_2^2}{1-\sigma_2^2},\quad
 q_X:=\frac{1+\rho^2}{2},\quad
 c_X:=\kappa_M^2\frac{1+\rho^2}{1-\rho^2}.
\]
In particular, $q_T,q_X<1$. The proof is given in
Appendix~\ref{appen}.
\end{lemma}

\begin{proof}[Proof of Theorem~\ref{thm:dprgt}]
For brevity set
\[
 a_X=64\bar L^2,\quad a_S=16\bar L^2,\quad
 b_E=3h_T,\quad b_X=96h_T\bar L^2,\quad
 b_S=24h_T\bar L^2,
\]
and let $\delta_T=1-q_T$ and $\delta_X=1-q_X$. Choose constants in the
following order:
\begin{equation}\label{eq:potential-weights}
 A=\frac72,\qquad B_0=\frac{5}{2\delta_T},\qquad
 C_0=\frac{Aa_X+2}{\delta_X},\qquad
 \kappa\ge16(Aa_S+C_0c_X).
\end{equation}
After these choices, take the constant $\theta$ in
\eqref{eq:repaired-parameters} sufficiently small so that all local
projection, descent, and coefficient restrictions below hold. Define
\begin{equation}\label{eq:coupled-potential}
 \Phi_k:=\mathbb{E} f(\bar x_k)
 +A\frac\alpha\tau\mathsf{E}_k+B_0\alpha\mathsf{T}_k
 +C_0\frac\alpha\tau\mathsf{X}_k.
\end{equation}
Using Lemma~\ref{lem:coupled-one-step}, $\tau=\kappa\alpha^2$, and
$1-(1-\tau)^2=\tau(2-\tau)$, the coefficients of
$\mathsf{E}_k,\mathsf{T}_k,\mathsf{X}_k$, and $\mathsf{S}_k$ in
$\Phi_{k+1}-\Phi_k$ are bounded above by
\begin{align*}
 &\alpha\{\tfrac32-A(2-\tau)+B_0b_E\tau^2\},
 &&-\alpha\{B_0\delta_T-\tfrac32\},\\
 &\frac\alpha\tau\{-C_0\delta_X+Aa_X
       +\tau(d_X+B_0b_X)\},
 &&\alpha\{-\tfrac14+(Aa_S+C_0c_X)/\kappa+B_0b_S\alpha^2\}.
\end{align*}
Thus, for this sufficiently small $\theta$, there is a constant $C_\Phi>0$
such that
\begin{equation}\label{eq:master-drift}
 \Phi_{k+1}-\Phi_k
 \le-\frac\alpha4\mathsf{G}_k-\frac\alpha8\mathsf{S}_k
     -\alpha\mathsf{E}_k-\alpha\mathsf{T}_k
     -\frac\alpha\tau\mathsf{X}_k+C_\Phi\alpha^3.
\end{equation}
The last term contains the oracle contributions because
$\alpha\tau\nu^2=\kappa\nu^2\alpha^3$; all constants are independent of
$k$ and $K$.

It remains to handle the one-sample initialization. Let
$N=\lfloor K/2\rfloor$. The variance assumption gives
$\mathsf{E}_0\le\nu^2$, consensus initialization gives $\mathsf{X}_0=0$, and
$\bs_0=\bd_0$ together with $\|d_{i,0}\|\le B$ gives
$\mathsf{T}_0=n^{-1}\mathbb E
\|((I-J)\otimes I_d)\bd_0\|^2\le B^2$.
Consequently, summing \eqref{eq:master-drift} from $0$ to $N$ gives
\[
 \alpha\sum_{k=0}^{N}\mathsf{E}_k
 \le\Phi_0-f_*+C_\Phi\alpha^3(N+1)=\mathcal O(\alpha^{-1}).
\]
Here $\alpha^3K=\mathcal O(1)$. Since
$N+1=\Theta(\alpha^{-3})$, there is a deterministic index $j\le N$ such
that $\mathsf{E}_j=\mathcal O(\alpha)$. Compactness, the uniform bounds on
$\mathsf{S}_j$ and $\mathsf{T}_j$, and Theorem~\ref{lem:consensus} then give
\[
 \Phi_j-f_*=\mathcal O(1):
 \quad \frac\alpha\tau\mathsf{E}_j=\mathcal O(1),\qquad
 \alpha\mathsf{T}_j=\mathcal O(\alpha),\qquad
 \frac\alpha\tau\mathsf{X}_j=\mathcal O(\alpha).
\]
Summing \eqref{eq:master-drift} again from $j$ to $K$ yields
\[
 \frac\alpha4\sum_{k=j}^{K}\mathsf{G}_k
 \le\Phi_j-f_*+C_\Phi\alpha^3(K-j+1)=\mathcal O(1).
\]
Since $K-j+1\ge K/2$, we obtain
\[
 \min_{j\le k\le K}\mathsf{G}_k
 =\mathcal O((\alpha K)^{-1})=\mathcal O(K^{-2/3}),
\]
which proves \eqref{eq:repaired-optimality-rate}.
\end{proof}

\begin{proof}[Proof of Corollary~\ref{coll}]
Choose an index $k_*$ attaining the minimum in
\eqref{eq:repaired-optimality-rate}. Theorem~\ref{lem:consensus} gives the
consensus estimate at the same index. Since Algorithm~\ref{alg:drgta} uses
two same-sample stochastic-gradient evaluations per node in each recursive
update (and one at initialization), its per-node SFO cost is $\mathcal O(K)$.
Taking $K=\mathcal O(\epsilon^{-3/2})$ completes the proof.
\end{proof}
\end{revisionblock}

\section{Numerical experiments}\label{sec:exm}
In this section, we compare our proposed DPRSRM with DRSGD in \cite{chen2021decentralized} and DRPGD in \cite{deng2023decentralized} on decentralized principal component analysis. It is important to note that the original DRPGD is a deterministic algorithm that utilizes the local full gradient. To adapt it to the online setting, we replace the full gradient with a stochastic gradient for the local updates. The numerical results on decentralized low-rank matrix completion are provided in the  supplementary 
 material.

\subsection{Decentralized principal component analysis} \label{sub:pca}
The decentralized principal component analysis (PCA) problem can be expressed mathematically as follows:
\be \label{prob:pca}
\min _{\mathbf{x} \in \mathcal{M}^{n}}-\frac{1}{2 n} \sum_{i=1}^{n} \text{tr}(x_{i}^{\top} A_{i}^{\top} A_{i} x_{i}), \;\; \text { s.t. } \;\; x_{1}=\ldots=x_{n},
\ee
where \(\mathcal{M}\) is the Stiefel manifold \(\text{St}(d,r)\), \(A_{i} \in \mathbb{R}^{m_{i} \times d}\) is the data matrix corresponding to the \(i\)-th node, and \(m_{i}\) denotes the number of samples. It is worth noting that if \(x^*\) is a solution to this problem, then any transformation of \(x^*\) by an orthogonal matrix \(Q \in \mathbb{R}^{r \times r}\) is also a valid solution. The distance between two points \(x\) and \(x^*\) is then calculated as:

\[
d_s(x, x^*) := \min_{Q \in \mathbb{R}^{r \times r}, \; Q^\top Q = QQ^\top = I_r} \; \|xQ - x^*\|.
\]

\subsubsection{Synthetic dataset}
In our study, we set the parameters as follows: \(m_{1} = \ldots = m_{n} = 1000\), \(d = 10\), and \(r = 5\). A matrix \(B \in \mathbb{R}^{1000n \times d}\) is generated, and its singular value decomposition (SVD) is performed, yielding \(B = U \Sigma V^\top\), where \(U \in \mathbb{R}^{1000n \times d}\) and \(V \in \mathbb{R}^{d \times d}\) are orthogonal matrices, and \(\Sigma \in \mathbb{R}^{d \times d}\) is a diagonal matrix. To control the distribution of singular values, we define \(\tilde{\Sigma} = \text{diag}(\gamma^j)\) with \(\gamma\) chosen from the interval \((0,1)\). The matrix \(A\) is then formed as \(A = U \tilde{\Sigma} V^\top \in \mathbb{R}^{1000n \times d}\). The matrices \(A_i\) are derived by partitioning the rows of \(A\) into \(n\) equally sized subsets. It can be shown that the first \(r\) columns of \(V\) represent the solution to \eqref{prob:pca}. In our experiments, the parameters \(\gamma\) and \(n\) are set to 0.8 and 8, respectively.

We employ fixed step sizes for all algorithms. The step size is set to $\alpha = \frac{\hat{\beta}}{\sqrt{K}}$ with $K$ being the maximal number of iterations.  The grid search is utilized to find the best $\hat{\beta}$ for each algorithm. The momentum parameter is chosen as $\tau = 0.999$. The batch size in each node is set as $10$ and the maximum iteration is set $K = 2000$. The clipping constant is set as $B=10^{8}$.    We choose the polar decomposition as the retraction operator for DRSGD. We test several graph matrices to model the topology across the nodes, namely, the Erdos-Renyi (ER) network with probability $p = 0.3, 0.6$, and the Ring network. Throughout this section, we select the mixing matrix $W$ to be the Metropolis constant edge weight matrix \cite{shi2015extra}.

\begin{figure}[htp]
	\centering
	\includegraphics[width = 0.46 \textwidth]{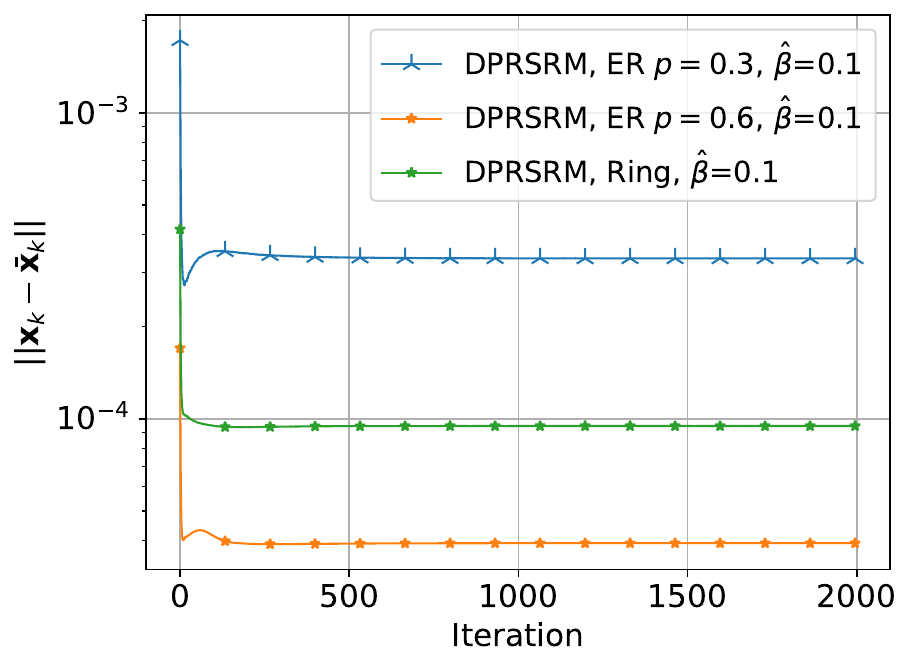}
	\includegraphics[width = 0.46 \textwidth]{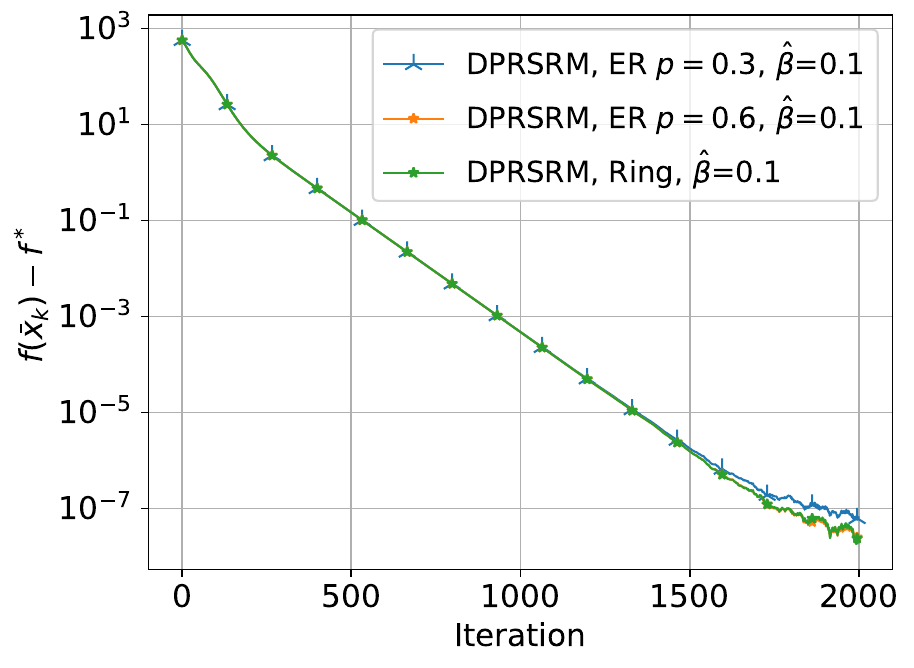} \\
	\includegraphics[width = 0.46 \textwidth]{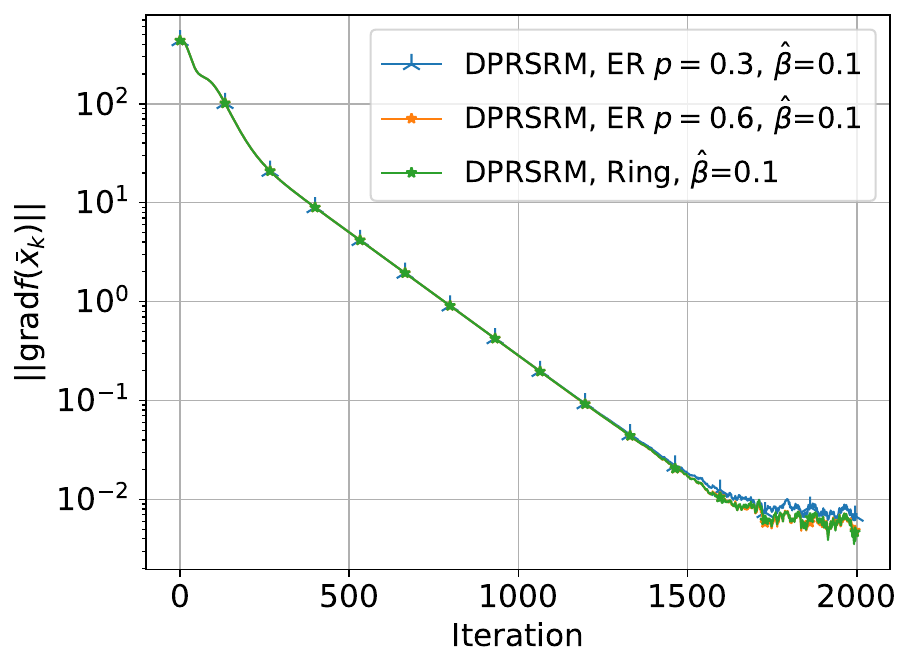}
	\includegraphics[width = 0.46 \textwidth]{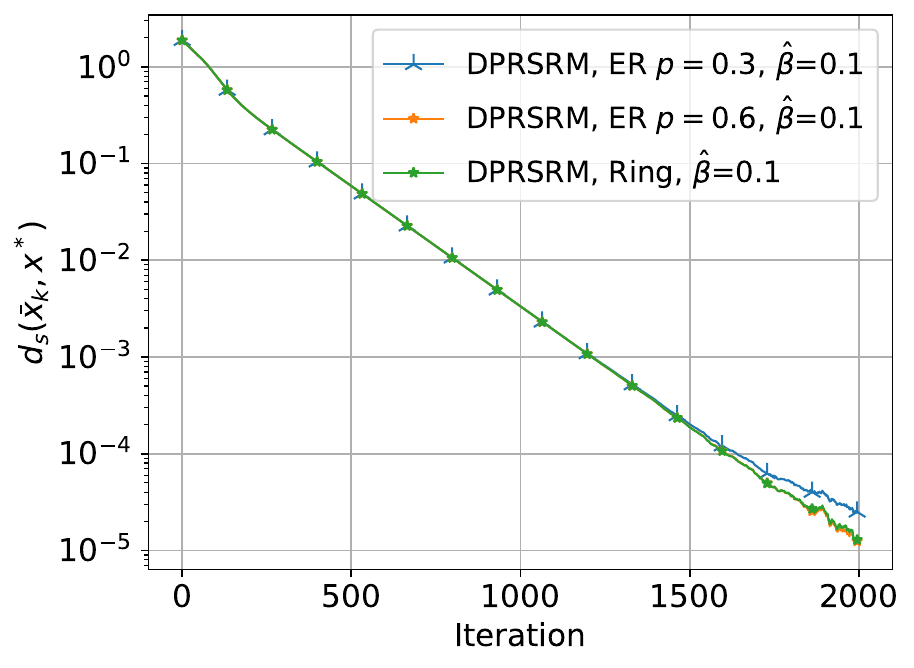}
	\caption{Numerical results on the synthetic dataset with different network graphs.}	
	\label{fig:num-pca-graph}
\end{figure}

\begin{figure}[htp]
	\centering
    \includegraphics[width = 0.46 \textwidth]{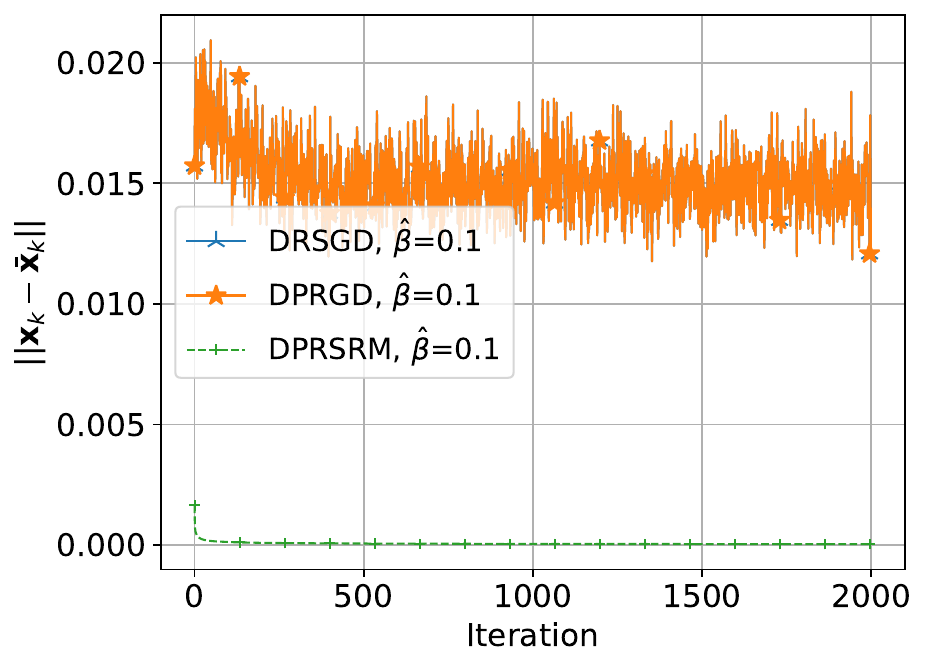}
    \includegraphics[width = 0.46 \textwidth]{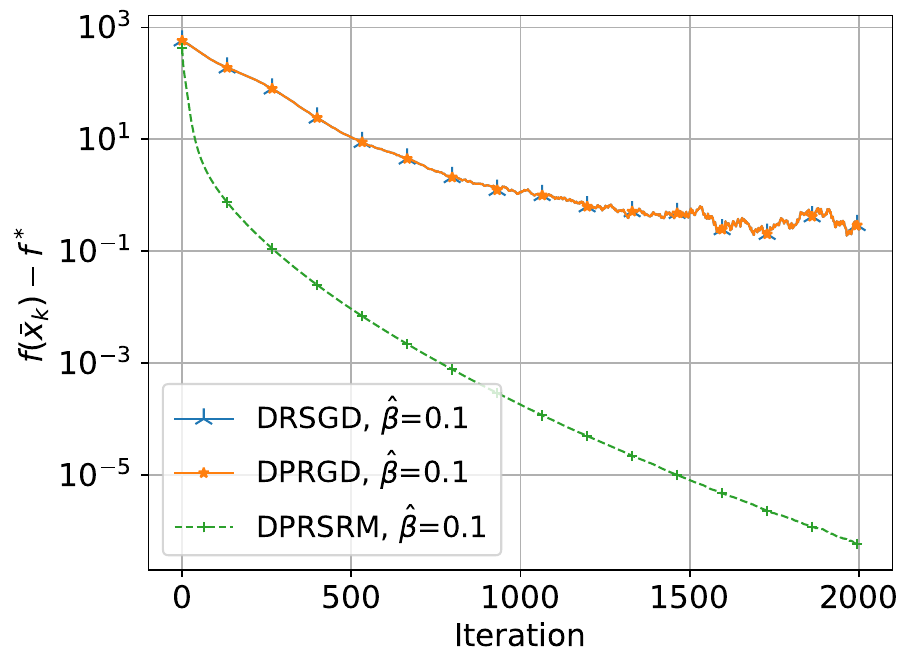} \\
    \includegraphics[width = 0.46 \textwidth]{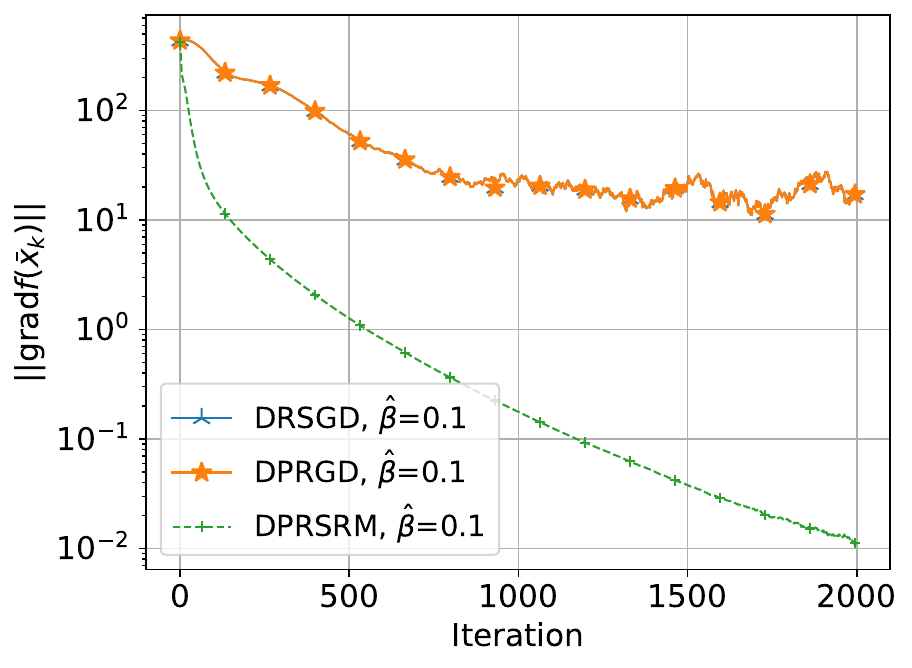}
    \includegraphics[width = 0.46 \textwidth]{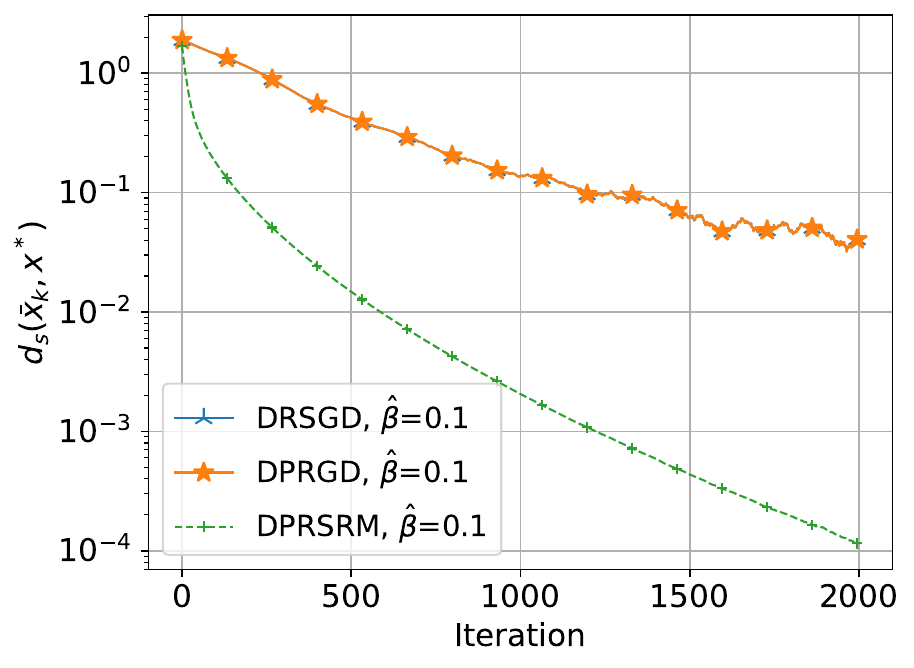}
	\caption{Results on the synthetic dataset with ER $p=0.6$.}	
	\label{fig:num-pca-alg}
\end{figure}

Firstly, we test all algorithms with different network graphs, namely, Ring, ER $p=0.3$, and ER $p=0.6$. The results are shown in Figure \ref{fig:num-pca-graph}. We see that there is not much difference among different graphs except for the consensus error. This is consistent with the existing results for DRSGD \cite{chen2021decentralized} and DPRGD \cite{deng2023decentralized}. 
Secondly, Figure \ref{fig:num-pca-alg} presents a comparison among the three algorithms. Our DPRSRM outperforms the other two, with DRSGD and DPRGD showing comparable performance.


\begin{figure}[htp]
	\centering
	\includegraphics[width = 0.46 \textwidth]{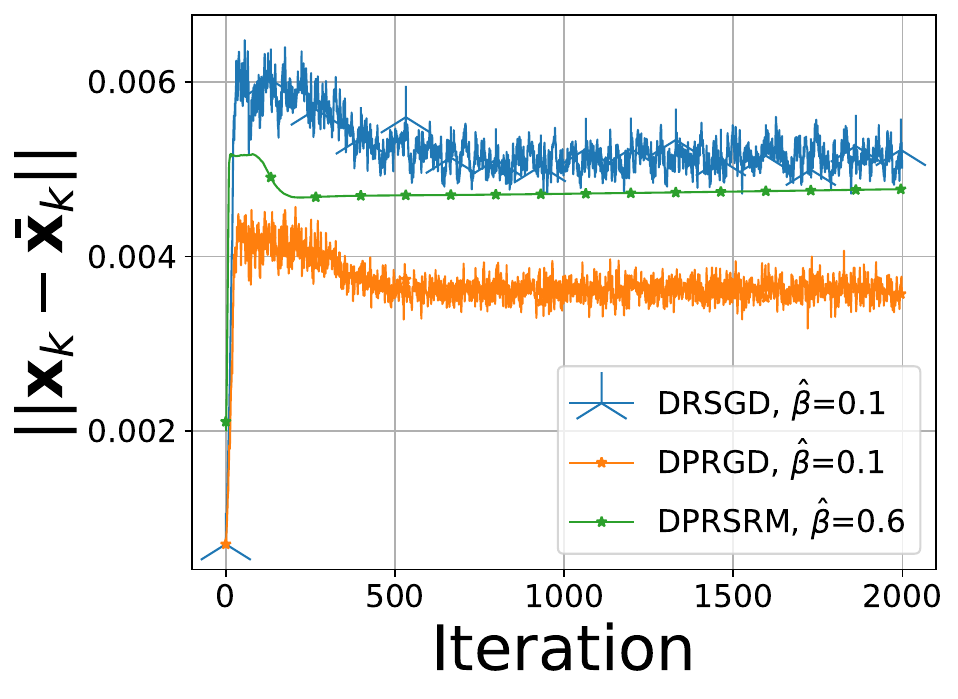}
	\includegraphics[width = 0.46 \textwidth]{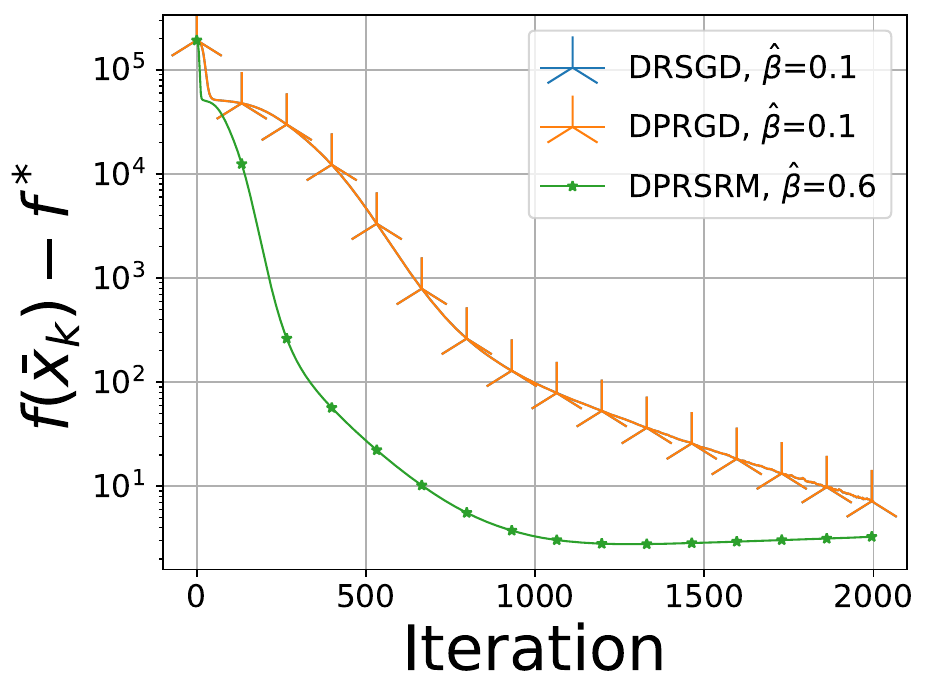}\\
    \includegraphics[width = 0.46 \textwidth]{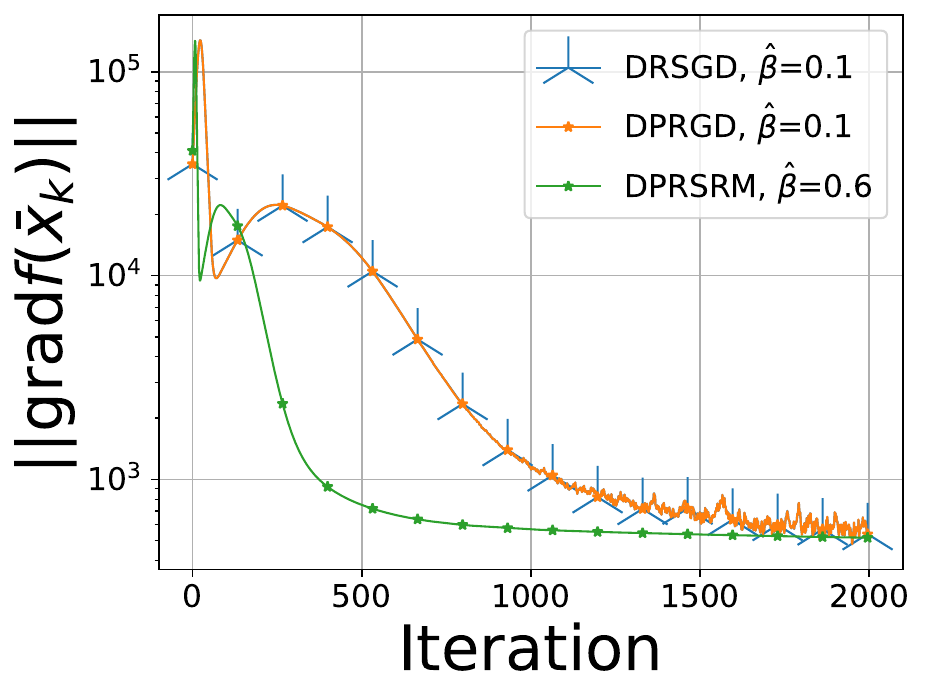}
	\includegraphics[width = 0.46 \textwidth]{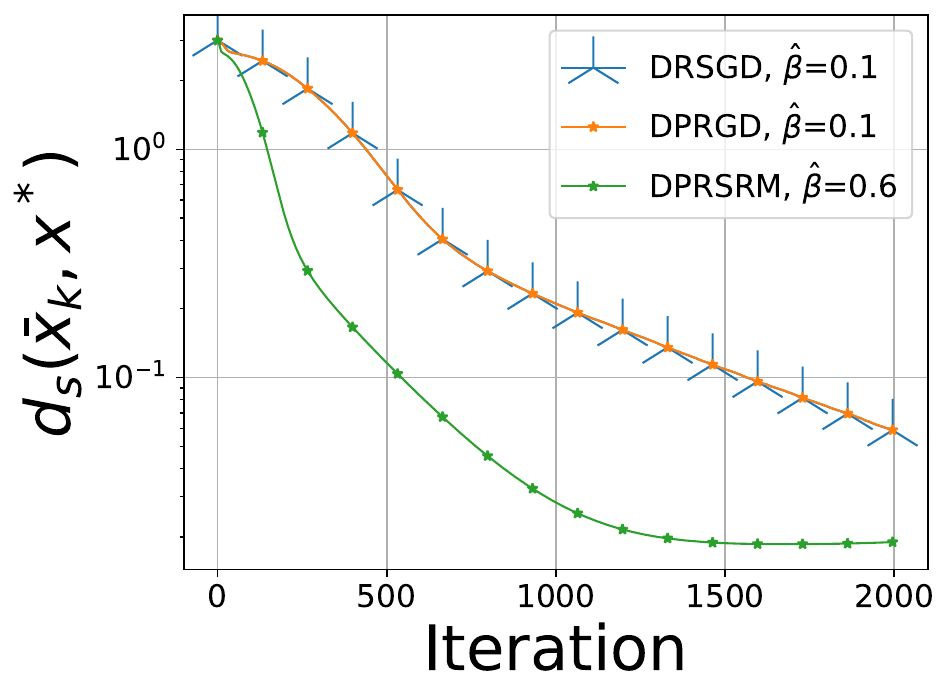}\\
	\caption{Results on Mnist dataset with ER $p=0.3$.}	
	\label{fig:num-pca-real}
\end{figure}


\subsubsection{Mnist dataset}
We also conduct numerical tests on the Mnist dataset \cite{lecun1998mnist}. The training images consist of 60000 handwritten images of size $32 \times 32$ and are used to generate $A_i$'s. We first normalize the data matrix by dividing 255 and randomly split the data into $n=8$ nodes with equal cardinality. Then, each node holds a local matrix $A_i$ of dimension $\frac{60000}{n} \times 784$. We compute the first 5 principal components, i.e., $d=784, r=5$.

For all algorithms, we use the fixed step sizes $\alpha = \frac{\hat{\beta}}{60000}$ with a best-chosen $\hat{\beta}$, batch size 1500 and momemtum parameter $\tau = 0.999$. 
Similar to the synthetic setting, our DPRSRM algorithm demonstrates superior performance compared to other algorithms in terms of objective function values, gradient norms, and distances to the optimal solution. It is important to note that due to the use of stochastic gradients, consensus among the algorithms may be affected by noise, which can be reduced by using a smaller step size.




\subsection{Low-rank matrix completion problem}\label{sec:more}
The low-rank matrix completion (LRMC) problem aims to reconstruct a matrix \(A \in \mathbb{R}^{d \times T}\) with low rank from its partially observed entries. Let \(\Omega\) represent the set of indices corresponding to the observed entries in \(A\). The LRMC problem of rank \(r\) can be expressed as:
\be \label{prob:lrmc} \min_{X \in {\rm Gr}(d,r), V\in \R^{r \times T}} \quad \frac{1}{2} \| \Pcal_{\Omega}(XV- A) \|^2, \ee
where \(\text{Gr}(d,r)\) represents the Grassmann manifold of \(r\)-dimensional subspaces in \(\mathbb{R}^d\), and \(\mathcal{P}_{\Omega}\) is a projection operator that selects the elements of \(A\) indexed by \(\Omega\), setting the rest to zero.

In a decentralized context, assume the matrix \(\mathcal{P}_{\Omega}(A)\) is divided into \(n\) equal parts by columns, labeled as \(A_1, A_2, \ldots, A_n\), each part corresponding to a different node. Replacing the Grassmann manifold constraint with the Stiefel manifold, the decentralized LRMC problem \cite{deng2023decentralized} is therefore formulated as:
\be \label{prob:dlrmc}
\begin{aligned}
\min \quad & \frac{1}{2} \sum_{i=1}^n \| \Pcal_{\Omega_i}(X_i V_i(X) - A_i) \|^2,  \\
\st \quad & X_1 = X_2 = \cdots = X_n, \\
& X_i \in {\rm St}(d,r), \quad \forall i \in [n],
\end{aligned}
\ee
where $\Omega_i$ is the corresponding indices set of $\Omega$ and $V_i(X):= {\rm argmin}_{V} \| \Pcal_{\Omega_i} (XV - A_i)\|$.

\begin{figure}[htp]
	\centering
	\includegraphics[width = 0.32\textwidth]{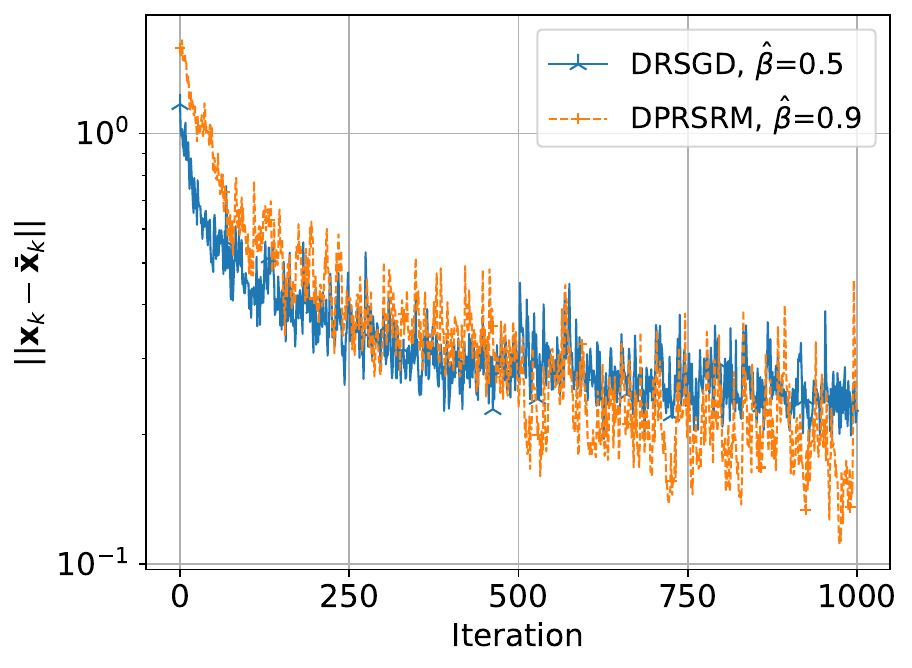}
    \includegraphics[width = 0.32\textwidth]{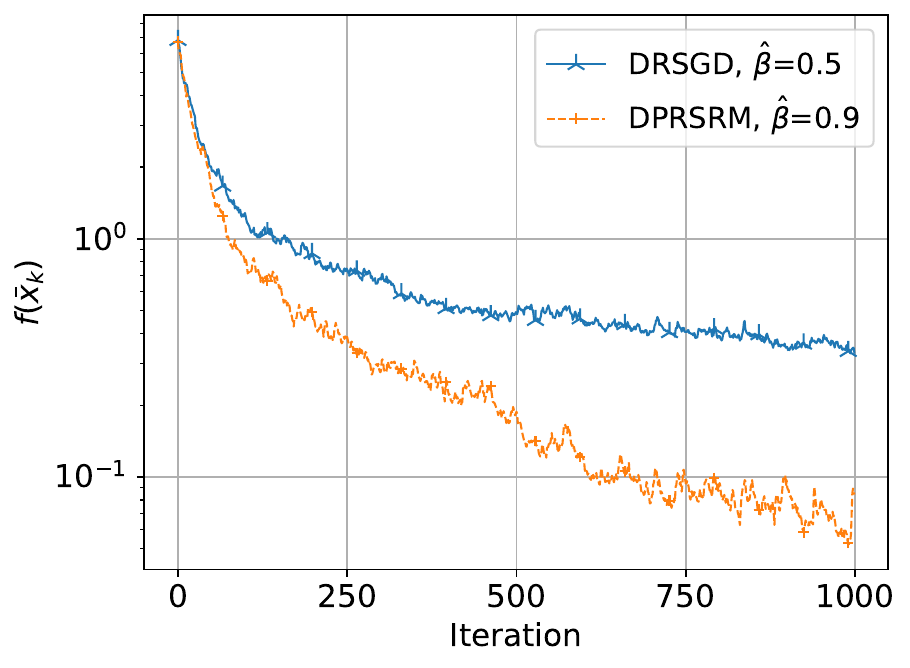}
     \includegraphics[width = 0.32\textwidth]{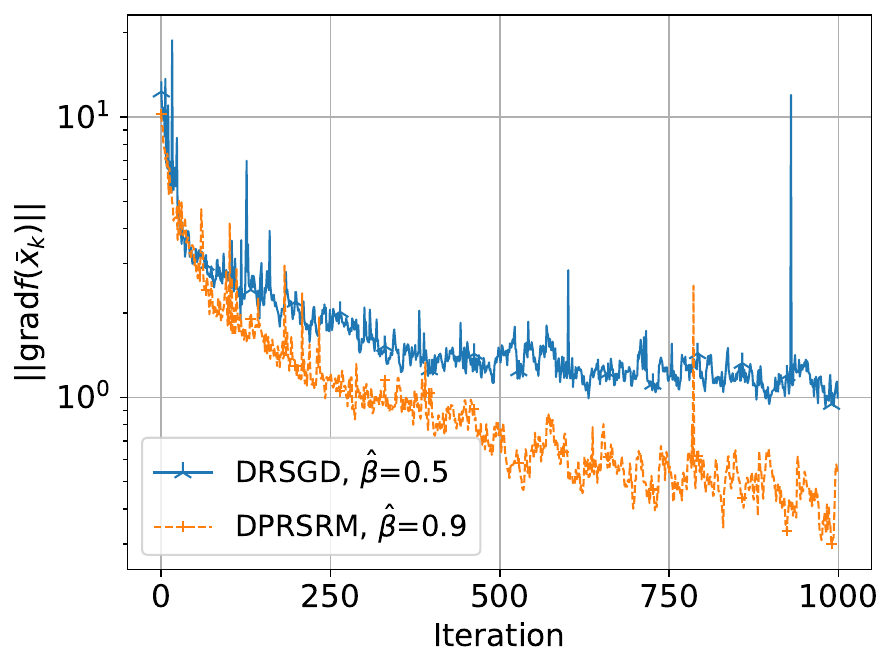}
	\caption{Numerical results for the decentralized LRMC problem with the Ring graph.}
	\label{fig:num-mc}	
\end{figure}

For numerical tests, we consider random generated $A$. To be specific, we first generate two random matrices $L \in \R^{d \times r}$ and $R \in \R^{r \times T}$, where each element obeys the standard Gaussian distribution. For the indices set $\Omega$, we generate a random matrix $B$ with each element following from the uniform distribution, then set $\Omega_{ij} = 1$ if $B_{ij} \leq \nu$ and $0$ otherwise. The parameter $\nu$ is set to $r(d+T-r)/(dT)$. In the implementations, we set $T=1000, d=50, r=10$, and $\alpha = \frac{\hat{\beta}}{\sqrt{K}}$ for all algorithms with $K$ being the maximal number of iterations. $\hat{\beta}$ is tuned to get the best performance for each algorithm individually. The Ring graph is used. The results are reported in Figure \ref{fig:num-mc}, where DPRGD is omitted due to its similar performance with DRSGD. We see that DPRSRM outperforms DRSGD.

\section{Conclusions and Limitations}
This work develops a decentralized projection Riemannian stochastic recursive momentum method by assuming that each
node has access to a stochastic first-order oracle.  Our algorithm
leverages local hybrid variance reduction and gradient tracking
to achieve a lower oracle complexity compared with the existing online
methods. It requires only $\mathcal{O}(1)$ gradient evaluations per iteration for each local node and does not require restarting with a large batch gradient.   


\section*{Acknowledgements} 

This effort was supported by the SciAI Center, and funded by the Office of Naval Research (ONR), under Grant Number N00014-23-
1-2729 (J.H.), as well as by the National Natural Science Foundation of China
(NSFC) grants 12401419 (K.D.).

\bibliographystyle{siamplain}
\bibliography{ref}

\newpage
\appendix
\section{Technical Lemmas}\label{appen}

\begin{lemma}\label{lem:inequa}
    Given any vectors $a,b\in \mathbb{R}^n$, it holds that
    \begin{equation}
        \langle a, \frac{b}{\|b\|} \rangle \geq \|a\| - 2\|a - b\|.
    \end{equation}
\end{lemma}

\begin{proof}
  It follows from the Cauchy inequality that
\begin{equation}
\begin{aligned}
    \langle a, \frac{b}{\|b\|} \rangle & =   \langle a - b, \frac{b}{\|b\|} \rangle + \|b\| \geq -\|a - b\| + \|b\| \geq  - 2\| a-b \| + \|a\|, \\
    \end{aligned}
\end{equation}
where the second inequality use $\|a\| \leq \|a - b\| + \|b\|$.
\end{proof}

\begin{lemma}[\cite{xu2015augmented}, Lemma 2]\label{lem:curr}

Let $u_k$ and $w_k$ be two positive scalar sequences such that for all $k\geq 1$
\begin{equation}\label{2}
  u_{k} \leq \eta u_{k-1} + w_{k-1},
\end{equation}
where $\eta\in (0,1)$ is the decaying factor.  Then we have 
\begin{equation}\label{22}
       \sum_{k=0}^K u_k \leq \frac{u_0}{1-\eta} +  \frac{1}{1-\eta}  \sum_{k=0}^{K-1} w_{k}. 
\end{equation}
\end{lemma}


The following inequality is the control of the distance between the Euclidean mean $\hat{x}$ and the manifold mean $\bar{x}$ by the square of consensus error.
\begin{lemma}[\cite{deng2023decentralized}]\label{lemma:quadratic}
   For any $\bx\in \Mcal^n$ satisfying $\|x_i - \bar{x}\| \leq {\delta}$, {$i\in [n]$,}  we have 
   \be\label{eq:distance-an-rm}
   \|\bar{x} - \hat{x} \| \leq M_2 \frac{\|\bx - \bar{\bx}\|^2}{n},
   \ee
   where $M_2 =\max_{x\in { \bar{U}_\Mcal({\delta})}} \|D^2 \Pcal_{\Mcal}(x) \|_{\rm op}$.
\end{lemma}

The following Lipschitz-type inequality for the projection operator $\Pcal_{\Mcal}(\cdot)$ is crucial in the analysis of projection-based methods. 

\begin{lemma}[\cite{deng2023decentralized}]\label{lemma-project}
For any $x\in\Mcal, u\in \{u\in \mathbb{R}^{d\times r}:\|u\|\leq {\delta}\}$, 
there exists a constant $Q$ such that
\be\label{projec-second-order1}
\| \Pcal_{\Mcal}(x + u)  - x - \Pcal_{T_x\Mcal}(u) \| \leq Q\|u\|^2.
\ee

\end{lemma}

\section{Proof of Section \ref{sec4}}

\subsection{Proof of Section \ref{sec:linear-consensus}}

\begin{revisionblock}
\begin{proof}[Proof of Lemma \ref{lemma:neibohood}]

The tracking identity and double stochasticity give
$\hat{\bs}_k=\hat{\bd}_k$ for every $k$. By the bounded initialization,
$\|\bs_0\|=\|\bd_0\|\leq\sqrt nB$. Suppose that
$\|\bs_k\|\leq3\sqrt nB/(1-\sigma_2)$.
Then, we have
\bee
\begin{aligned}
\| \bs_{k+1} - \hat{\bd}_k\| & = \| \bW \bs_k - \hat{\bd}_k + \bd_{k+1} - \bd_k \| \\
& =   \| \bW \bs_k - \hat{\bs}_k + \bd_{k+1} - \bd_k \|\\
& =  \| [(W-J)\otimes I_d ]  \bs_k  + \bd_{k+1} - \bd_k \|\\
& \leq  \sigma_2 \| \bs_k \| + 2 \sqrt{n} B.
\end{aligned}
\eee
Since $\|\hat{\bd}_k\|\leq\sqrt nB$, it follows that
\[
 \|\bs_{k+1}\|\leq\frac{3\sqrt nB\sigma_2}{1-\sigma_2}
 +3\sqrt nB=\frac{3\sqrt nB}{1-\sigma_2}.
\]
Induction proves \eqref{eq:sk-bound}.
\end{proof}

\begin{proof}[Proof of Lemma \ref{lem:stay-neighborhood}]
    We prove the two claims by induction. Assume that
    $\bx_k\in\mathcal{N}(\delta)$. For any $i\in[n]$,
\[ \begin{aligned}
\| \sum_{j=1}^n W_{ij} x_{j,k} - \alpha v_{i,k} - \bar{x}_k  \| 
 \leq &   \| \sum_{j=1}^n W_{ij} (x_{j,k} - \bar{x}_k ) - \alpha  v_{i,k}   \| \\
 \leq &   \sum_{j=1}^n W_{ij} \|x_{j,k} - \bar{x}_k  \| + \alpha \|v_{i,k} \| \\
 \leq   &  \delta +   \sqrt{nD_B}\alpha \leq 2\delta.
\end{aligned}\]
By $\delta<R/2$, we have
$\sum_{j=1}^nW_{ij}x_{j,k}-\alpha v_{i,k}\in\bar U_{\Mcal}(2\delta)$.
Moreover, $\max_i\|x_{i,k}-\bar x_k\|\leq\delta$, so
$\hat x_k\in\bar U_{\Mcal}(\delta)\subset\bar U_{\Mcal}(2\delta)$.
The $\kappa_M$-Lipschitz property of $\Pcal_{\Mcal}$ therefore gives
    \[ \begin{aligned}
      &  \| \bx_{k+1} - \bar{\bx}_{k+1} \|  \leq \| \bx_{k+1} - \bar{\bx}_k \| \\
        &= \|\Pcal_{\Mcal^n}( \bW \bx_k - \alpha \bv_k) - \Pcal_{\Mcal^n}(\hat{\bx}_k) \| \\
        & \leq \frac{R}{R- 2\delta} \| \bW \bx_k - \alpha \bv_k - \hat{\bx}_k\| \\
        & \leq \frac{R}{R- 2\delta} \| [ (W - J) \otimes I_d] ( \bx_k -\hat{\bx}_k) -\alpha \bv_k  \| \\
        & \leq \frac{R \sigma_2}{R- 2\delta} \|\bx_k - \bar{\bx}_k\| + \frac{R \alpha}{R- 2\delta} \|\bv_k\| \\
        & \leq \frac{R \sigma_2}{R- 2\delta} \delta + \frac{R \alpha}{R- 2\delta} \sqrt{nD_B}.
    \end{aligned} \]
 For a fixed $\delta$, the inequality
 $\|\bx_{k+1}-\bar{\bx}_{k+1}\|\leq\delta$ holds provided that
 \begin{equation}
     \alpha \leq \frac{(R(1-\sigma_2) - 2\delta)\delta}{R\sqrt{nD_B}}.
 \end{equation}
 The right-hand side is positive whenever
 $0<\delta<R(1-\sigma_2)/2$, as required by \eqref{cond-alpha-delta}.
 This closes the induction.
\end{proof}
\end{revisionblock}

\subsection{Proof of Section \ref{sec:proof:dprgt}}

\begin{revisionblock}
\begin{proof}[Proof of Lemma~\ref{lem:coupled-one-step}]
We prove the four inequalities in turn. For the transition from $k$ to
$k+1$, all conditional expectations are taken with respect to
$\mathcal F_{k+1}$, i.e., immediately before drawing
$\boldsymbol\xi_{k+1}$. Thus, $\bx_k,\bx_{k+1}$ and $\bd_k$ are measurable
with respect to the conditioning sigma-algebra.

\emph{Step 1: estimator error.}
Set
\[
 \Delta_k:=\bd_k-\grad f(\bx_k),\qquad
 \mathbf g_k^\xi:=\grad f(\bx_k,\boldsymbol\xi_{k+1}),\qquad
 \mathbf g_k:=\grad f(\bx_k).
\]
Since $B\geq G_{\max}$, every block of $\mathbf g_{k+1}$ belongs to the clipping
ball. Blockwise nonexpansiveness of the Euclidean projection therefore gives
\begin{equation}\label{eq:clip-estimator-comparison}
 \|\bd_{k+1}-\mathbf g_{k+1}\|^2
 \leq\|\bq_{k+1}-\mathbf g_{k+1}\|^2.
\end{equation}
The un-clipped error has the exact decomposition
\begin{equation}\label{eq:estimator-centered-decomp}
 \bq_{k+1}-\mathbf g_{k+1}=(1-\tau)\Delta_k
 +\boldsymbol\varepsilon_{k+1},
\end{equation}
where
\begin{align*}
 \boldsymbol\varepsilon_{k+1}:={}&
 \tau(\mathbf g_{k+1}^\xi-\mathbf g_{k+1})\\
 &+(1-\tau)\{\mathbf g_{k+1}^\xi-\mathbf g_k^\xi
                  -(\mathbf g_{k+1}-\mathbf g_k)\}.
\end{align*}
Both terms in braces are evaluated using the same fresh sample. Their sum is
conditionally centered, so
$\mathbb E[\boldsymbol\varepsilon_{k+1}\mid\mathcal F_{k+1}]=0$ and the cross term with
$\Delta_k$ in \eqref{eq:estimator-centered-decomp} vanishes. We do not assume
that the two fresh-sample terms are mutually independent. Instead,
$\|a+b\|^2\leq2\|a\|^2+2\|b\|^2$, conditional variance is bounded by the
corresponding second moment, and Assumption~\ref{assum:stochatis} give
\begin{align}
 \mathbb E[\|\boldsymbol\varepsilon_{k+1}\|^2\mid\mathcal F_{k+1}]
 \leq{}&2n\tau^2\nu^2
 +2(1-\tau)^2\bar L^2\|\bx_{k+1}-\bx_k\|^2.       \label{eq:centered-noise-bound}
\end{align}
Taking full expectation in \eqref{eq:clip-estimator-comparison} and using
\eqref{eq:estimator-centered-decomp}--\eqref{eq:centered-noise-bound} yields
\begin{equation}\label{eq:repaired-estimator-step}
 \mathsf{E}_{k+1}\leq(1-\tau)^2\mathsf{E}_k+2\tau^2\nu^2
 +\frac{2\bar L^2}{n}\mathbb{E}\|\bx_{k+1}-\bx_k\|^2.
\end{equation}
This is the point at which the original proof introduced the spurious factor
$\alpha^2$: mean-square Lipschitz continuity gives the last squared
iterate-difference directly.

For completeness, let $\by_k:=\bW\bx_k-\alpha\bv_k$. Since
$\bx_{k+1}=\Pcal_{\Mcal^n}(\by_k)$ and $\bx_k\in\Mcal^n$, optimality of the
projection gives $\|\bx_{k+1}-\by_k\|\leq\|\bx_k-\by_k\|$. Hence,
\begin{align*}
 \|\bx_{k+1}-\bx_k\|
 &\leq2\|\by_k-\bx_k\|\\
 &\leq2\|(I_{nd}-\bW)(\bx_k-\bar\bx_k)+\alpha\bv_k\|.
\end{align*}
Because $\|I_{nd}-\bW\|\leq2$ and
$\|\bv_k\|\leq\|\bs_k\|$, squaring the preceding display gives
\begin{equation}\label{eq:repaired-iterate-difference}
 \frac1n\mathbb{E}\|\bx_{k+1}-\bx_k\|^2
 \leq32\mathsf{X}_k+8\alpha^2\mathsf{S}_k.
\end{equation}
Combining \eqref{eq:repaired-estimator-step} and
\eqref{eq:repaired-iterate-difference}, and using $(1-\tau)^2\leq1$ in the
forcing term, proves \eqref{eq:coupled-E-main}.

\emph{Step 2: tracking error.}
The identity $d_{i,k+1}=q_{i,k+1}$ need not hold when clipping is active.
However, $d_{i,k}\in B\mathbb B$ for every $k\geq0$ (at $k=0$ this follows
from Assumption~\ref{assum:stochatis}), and therefore
\begin{align}
 \|d_{i,k+1}-d_{i,k}\|
 &=\|\Pcal_{B\mathbb B}(q_{i,k+1})-
       \Pcal_{B\mathbb B}(d_{i,k})\|\nonumber\\
 &\leq\|q_{i,k+1}-d_{i,k}\|.                 \label{eq:clip-increment}
\end{align}
Expanding the last vector gives three terms:
\begin{align*}
 q_{i,k+1}-d_{i,k}
 ={}&\grad f_i(x_{i,k+1},\xi_{i,k+1})
       -\grad f_i(x_{i,k},\xi_{i,k+1})\\
 &+\tau\{\grad f_i(x_{i,k},\xi_{i,k+1})-\grad f_i(x_{i,k})\}
   +\tau\{\grad f_i(x_{i,k})-d_{i,k}\}.
\end{align*}
Thus, by $\|a+b+c\|^2\leq3(\|a\|^2+\|b\|^2+\|c\|^2)$, the oracle
assumptions, and \eqref{eq:repaired-iterate-difference},
\begin{equation}\label{eq:repaired-difference-d}
 \frac1n\mathbb{E}\|\bd_{k+1}-\bd_k\|^2
 \leq96\bar L^2\mathsf{X}_k+24\bar L^2\alpha^2\mathsf{S}_k
 +3\tau^2\mathsf{E}_k+3\tau^2\nu^2.
\end{equation}

Double stochasticity and $\hat\bs_k=\hat\bd_k$ imply
\begin{equation}\label{eq:tracking-exact-decomp}
 \bs_{k+1}-\hat{\bd}_{k+1}
 =((W-J)\otimes I_d)(\bs_k-\hat{\bd}_k)
  +((I-J)\otimes I_d)(\bd_{k+1}-\bd_k).
\end{equation}
For $\sigma_2>0$, apply
$\|a+b\|^2\leq(1+\zeta)\|a\|^2+(1+1/\zeta)\|b\|^2$ with
\[
 \zeta_T:=\frac{1-\sigma_2^2}{2\sigma_2^2}.
\]
Then
\[
 (1+\zeta_T)\sigma_2^2=q_T,
 \qquad 1+\frac1{\zeta_T}=h_T.
\]
Using $\|I-J\|=1$ in \eqref{eq:tracking-exact-decomp} and then
\eqref{eq:repaired-difference-d} proves \eqref{eq:coupled-T-main}. If
$\sigma_2=0$, the first term in \eqref{eq:tracking-exact-decomp} vanishes;
the same displayed bound remains valid because $q_T=1/2$ and $h_T=1$.

\emph{Step 3: dynamic consensus error.}
The pathwise projection estimate already established in the proof of
Theorem~\ref{lem:consensus} gives
\begin{equation}\label{eq:pathwise-dynamic-consensus}
 \|\bx_{k+1}-\bar\bx_{k+1}\|
 \leq\rho\|\bx_k-\bar\bx_k\|+\kappa_M\alpha\|\bs_k\|.
\end{equation}
Minkowski's inequality in $L^2$, divided by $\sqrt n$, yields
\[
 \sqrt{\mathsf{X}_{k+1}}
 \leq\rho\sqrt{\mathsf{X}_k}+\kappa_M\alpha\sqrt{\mathsf{S}_k}.
\]
When $\rho>0$, Young's inequality with
$\zeta_X=(1-\rho^2)/(2\rho^2)$ gives
$(1+\zeta_X)\rho^2=q_X$ and
$(1+1/\zeta_X)\kappa_M^2=c_X$, proving
\eqref{eq:coupled-X-main}. When $\rho=0$, direct squaring gives
$\mathsf X_{k+1}\leq\kappa_M^2\alpha^2\mathsf S_k$, which is stronger than
\eqref{eq:coupled-X-main} with $q_X=1/2$ and $c_X=\kappa_M^2$.

\emph{Step 4: descent with explicit projection remainders.}
We now derive \eqref{eq:coupled-D-main} without hiding the geometric terms.
Define the ambient consensus residual and its tangent projection by
\[
 c_{i,k}:=x_{i,k}-\sum_{j=1}^nW_{ij}x_{j,k},\qquad
 \grad\phi_i(\bx_k):=\Pcal_{T_{x_{i,k}}\Mcal}c_{i,k}.
\]
Let $\mathbf c_k$ denote the stack of the $c_{i,k}$. Then
$\sum_i c_{i,k}=0$ and
$\|\mathbf c_k\|=\|((I-W)\otimes I_d)\bx_k\|
\leq2\|\bx_k-\bar\bx_k\|$. Since the tangent projector is Lipschitz on the
compact submanifold, there is $H_T>0$ such that
\begin{equation}\label{eq:tangent-projector-lipschitz}
 \|\Pcal_{T_x\Mcal}-\Pcal_{T_y\Mcal}\|_{\rm op}
 \leq H_T\|x-y\|,\qquad x,y\in\Mcal.
\end{equation}
Consequently,
\begin{align}
 \left\|\sum_{i=1}^n\grad\phi_i(\bx_k)\right\|
 &=\left\|\sum_{i=1}^n
   (\Pcal_{T_{x_{i,k}}\Mcal}-\Pcal_{T_{\bar x_k}\Mcal})c_{i,k}\right\|
 \nonumber\\
 &\leq2H_T\|\bx_k-\bar\bx_k\|^2.             \label{eq:sum-consensus-grad}
\end{align}

Put $w_{i,k}:=c_{i,k}+\alpha v_{i,k}$. By
\eqref{eq:xk-barxk-bound} and \eqref{eq:sk-bound},
\[
 \max_{1\leq i\leq n}\|w_{i,k}\|
 \leq C_w\alpha,
 \qquad C_w:=2\sqrt{n C_B}+\sqrt{n D_B}.
\]
After reducing $\theta$ if necessary, this quantity lies in the uniform
radius of the second-order projection estimate
\eqref{projec-second-order1}. Define
\[
 r_{i,k}:=x_{i,k+1}-x_{i,k}+\grad\phi_i(\bx_k)+\alpha v_{i,k}.
\]
Because $v_{i,k}\in T_{x_{i,k}}\Mcal$,
\begin{equation}\label{eq:projection-remainder-block}
 \|r_{i,k}\|\leq Q\|w_{i,k}\|^2.
\end{equation}
Writing $\hat r_k=n^{-1}\sum_i r_{i,k}$, using
\eqref{eq:sum-consensus-grad}, \eqref{eq:projection-remainder-block}, and
$\|\bv_k\|\leq\|\bs_k\|$, we obtain
\begin{align}
 \|\hat x_{k+1}-\hat x_k+\alpha\hat v_k\|
 &\leq\frac1n\left\|\sum_i\grad\phi_i(\bx_k)\right\|
       +\frac Qn\sum_i\|w_{i,k}\|^2\nonumber\\
 &\leq\frac{C_m}{n}\|\bx_k-\bar\bx_k\|^2
       +\frac{2Q\alpha^2}{n}\|\bs_k\|^2,       \label{eq:mean-step-remainder}
\end{align}
where $C_m:=2H_T+8Q$; here we used
$\sum_i\|w_{i,k}\|^2\leq8\|\bx_k-\bar\bx_k\|^2
+2\alpha^2\|\bs_k\|^2$.

Let
\[
 \delta_k:=\bar x_{k+1}-\bar x_k+\alpha\hat v_k,
 \qquad \mathcal C_k:=\|\bx_k-\bar\bx_k\|^2.
\]
The mean-versus-projected-mean estimate \eqref{eq:distance-an-rm} and
\eqref{eq:mean-step-remainder} give
\[
 \|\delta_k\|\leq
 \frac{M_2+C_m}{n}\mathcal C_k
 +\frac{M_2}{n}\mathcal C_{k+1}
 +\frac{2Q\alpha^2}{n}\|\bs_k\|^2.
\]
Therefore,
\begin{equation}\label{eq:delta-k-explicit}
 \|\delta_k\|^2\leq
 \frac{\Gamma_1}{n^2}\mathcal C_k^2
 +\frac{\Gamma_2}{n^2}\mathcal C_{k+1}^2
 +\frac{\Gamma_3\alpha^4}{n^2}\|\bs_k\|^4,
\end{equation}
where
\[
 \Gamma_1:=3(M_2+C_m)^2,\qquad
 \Gamma_2:=3M_2^2,\qquad \Gamma_3:=12Q^2.
\]
Moreover,
\begin{equation}\label{eq:bar-step-explicit}
 \|\bar x_{k+1}-\bar x_k\|^2
 \leq\frac{2\alpha^2}{n}\|\bs_k\|^2+2\|\delta_k\|^2.
\end{equation}

Let $g_k:=\grad f(\bar x_k)$ and
$\mathbf g_k^{\rm av}:=\mathbf1_n\otimes g_k$. The Riemannian quadratic
upper bound and the identity
\begin{equation}\label{eq:polarization-full-direction}
 -\alpha\langle g_k,\hat s_k\rangle
 =-\frac\alpha2\|g_k\|^2-\frac{\alpha}{2n}\|\bs_k\|^2
  +\frac{\alpha}{2n}\|\mathbf g_k^{\rm av}-\bs_k\|^2
\end{equation}
give
\begin{align}
 f(\bar x_{k+1})
 \leq{}&f(\bar x_k)-\frac\alpha4\|g_k\|^2
 -\frac{\alpha}{2n}\|\bs_k\|^2
 +\frac{\alpha}{2n}\|\mathbf g_k^{\rm av}-\bs_k\|^2\nonumber\\
 &+\frac1\alpha\|\delta_k\|^2
 +\frac L2\|\bar x_{k+1}-\bar x_k\|^2
 +\frac\alpha n\sum_{i=1}^n
   \langle g_k,s_{i,k}-v_{i,k}\rangle.             \label{eq:descent-before-bounds}
\end{align}
Indeed, \eqref{eq:descent-before-bounds} follows by writing
$\bar x_{k+1}-\bar x_k+\alpha\hat s_k
=\delta_k+\alpha(\hat s_k-\hat v_k)$ and applying Young's inequality only
to $\langle g_k,\delta_k\rangle$.

We next bound the two error terms in \eqref{eq:descent-before-bounds}. Define
$\hat g_k=n^{-1}\sum_i\grad f_i(x_{i,k})$ and
$\hat d_k=n^{-1}\sum_i d_{i,k}$. Jensen's inequality gives
\begin{align}
 \frac1n\mathbb E\|\mathbf g_k^{\rm av}-\bs_k\|^2
 &\leq3L_G^2\mathsf X_k+3\mathsf E_k+3\mathsf T_k.   \label{eq:average-gradient-error}
\end{align}
To see this, insert $\mathbf1_n\otimes\hat g_k$ and
$\hat\bd_k=\mathbf1_n\otimes\hat d_k$ between
$\mathbf g_k^{\rm av}$ and $\bs_k$; the three squared terms are bounded by
$L_G^2\mathsf X_k$, $\mathsf E_k$, and $\mathsf T_k$, respectively.

For the normal term, set
\[
 A_{i,k}:=g_k-\Pcal_{T_{x_{i,k}}\Mcal}g_k.
\]
Since $g_k\in T_{\bar x_k}\Mcal$, compactness and
\eqref{eq:tangent-projector-lipschitz} give
$\|A_{i,k}\|\leq H_N\|x_{i,k}-\bar x_k\|$ for a constant $H_N>0$.
Moreover, $s_{i,k}-v_{i,k}\in N_{x_{i,k}}\Mcal$ and
$\|s_{i,k}-v_{i,k}\|\leq\|s_{i,k}\|$. Thus
$\langle g_k,s_{i,k}-v_{i,k}\rangle
=\langle A_{i,k},s_{i,k}-v_{i,k}\rangle$, and blockwise Young's inequality
gives, for every fixed $\eta>0$,
\begin{equation}\label{eq:repaired-normal-term}
 \frac\alpha n\sum_{i=1}^n
 \mathbb E\langle g_k,s_{i,k}-v_{i,k}\rangle
 \leq\frac{\alpha\eta}{2}\mathsf S_k
 +\frac{\alpha H_N^2}{2\eta}\mathsf X_k.
\end{equation}
The factor $\alpha$ appears exactly once in this bound.

Finally, the pathwise estimates \eqref{eq:xk-barxk-bound} and
\eqref{eq:sk-bound} imply from \eqref{eq:delta-k-explicit} that
\begin{equation}\label{eq:delta-order-explicit}
 \mathbb E\|\delta_k\|^2\leq C_\delta\alpha^4,
 \quad
 C_\delta:= (\Gamma_1+\Gamma_2)C_B^2+\Gamma_3D_B^2.
\end{equation}
Choose $\eta=1/4$, $\alpha\leq1$, and $L\alpha\leq1/8$. Substituting
\eqref{eq:bar-step-explicit}, \eqref{eq:average-gradient-error},
\eqref{eq:repaired-normal-term}, and \eqref{eq:delta-order-explicit} into
\eqref{eq:descent-before-bounds} gives \eqref{eq:coupled-D-main} with, for
example,
\[
 d_X:=\frac32L_G^2+2H_N^2,
 \qquad C_D:=(1+L)C_\delta.
\]
All these constants depend only on the problem data and the fixed tube, and
are independent of $k,K,\alpha$, and $\tau$.
\end{proof}
\end{revisionblock}

\end{document}